\documentclass[12pt]{article}
\usepackage[english]{babel}
\usepackage[latin1]{inputenc}
\usepackage{amsfonts,amssymb,amsmath, epsfig}
\usepackage{color,graphicx,graphics,psfrag}
\usepackage{amsmath,amstext,amssymb,amsfonts, amscd}

\def\Box{\leavevmode\vbox{\hrule
     \hbox{\vrule\kern4pt\vbox{\kern4pt}%
           \vrule}\hrule}}
\def\blackbox{\leavevmode\vrule height 5pt width 4pt depth 0pt\relax}
\def\endproof{\null\hfill {$\blackbox$}\bigskip}

\newcounter{appendix}
\def\appendix{\advance\c@appendix by 1
   \def\thesection{\Alph{section}}
   \ifnum\c@appendix=1 \setcounter{section}{-1} \fi
   \@startsection {section}{1}{\z@}{-3.5ex plus -1ex minus 
   -.2ex}{2.3ex plus .2ex}{\Large\bf}}

\def\paragraph#1{{\bf #1\ }}

\newtheorem{lemma}{Lemma}[section]  

\newtheorem{theorem}[lemma]{Theorem}

\newtheorem{definition}[lemma]{Definition}

\newtheorem{proposition}[lemma]{Proposition}

\newtheorem{remark}{Remark}[section]

\newcommand{\dd}{\,\mathrm{d}}

\newcommand {\e} {\varepsilon}
\newcommand{\sgn}{\text{Sign}}

\title{A continuum model for nematic alignment of self-propelled particles} 
\author{Pierre Degond$^1$, Angelika Manhart$^2$, Hui Yu$^{1,3}$} 
\date{} 
\begin{document}

\maketitle

\begin{center}
1. Department of Mathematics, Imperial College London\\
London, SW7 2AZ, United Kingdom\\
pdegond@imperial.ac.uk 
\end{center}

\begin{center}
2. Faculty of Mathematics, University of Vienna, \\
Oskar-Morgenstern Platz 1,1090 Vienna, Austria,  \\
angelika.manhart@univie.ac.at
\end{center}

\begin{center}
3. Universit\'{e} de Toulouse; UPS, INSA, UT1, UTM\\
Institut de Math\'{e}matiques de Toulouse, France \\
and CNRS; Institut de Math\'{e}matiques de Toulouse, UMR 5219, France \\
hyu@math.univ-toulouse.fr
\end{center}

\vspace{0.5 cm}
\begin{abstract}
A continuum model for a population of self-propelled particles interacting through nematic alignment is derived from an individual-based model. The methodology consists of introducing a hydrodynamic scaling of the corresponding mean-field kinetic equation. The resulting perturbation problem is solved thanks to the concept of generalized collision invariants. It yields a hyperbolic but non-conservative system of equations for the  nematic mean direction of the flow and the densities of particles flowing parallel or anti-parallel to this mean direction. Diffusive terms are introduced under a weakly non-local interaction assumption and the diffusion coefficient is proven to be positive. An application to the modeling of myxobacteria is outlined. 
\end{abstract}

\medskip
\noindent
{\bf Acknowledgements:} This work has been supported by the Agence Nationale pour la Recherche (ANR) under grant MOTIMO (ANR-11-MONU-009-01), by the National Science Foundation (NSF) under grant RNMS11-07444 (KI-Net), by the Engineering and Physical Sciences Research Council (EPSRC) under grant ref: EP/M006883/1, by the Austrian Science Fund (FWF) through the doctoral school ``Dissipation and Dispersion in Nonlinear PDEs'' (project W1245) as well as the Vienna Science and Technology Fund (WWTF) (project LS13/029). PD is on leave from CNRS, Institut de Math\'ematiques de Toulouse, France. PD acknowledges support from the Royal Society and the Wolfson foundation through a Royal Society Wolfson Research Merit Award. 


\vskip 0.4cm

\setcounter{equation}{0}
\section{Introduction}
\label{sec:intro}

Systems of large numbers of self-propelled particles are commonly found in nature. In three space dimensions, well-known examples include fish schools, flocks of birds and insects \cite{Buhl2006, Cavagna2010, Parrish1997}. In two space dimensions, typical examples are bacteria, such as myxobacteria and \emph{Bacillus subtilis} \cite{Ben-Jacob2000, Kaiser2003, Sokolov2007, Welch2001, Zhang2010}. However, the interest in such systems is not limited to biological examples. It has been found that the underlying principles also often apply to human walking behavior \cite{Helbing2007, Moussaid2010} and physical systems, such as granular media \cite{Dhar2006, Narayan2007}. The main interest in these systems, both mathematically and from an application's point of view, stems from their emergent properties. Emergence refers to  the fact that fairly local interaction rules between the agents like attraction, alignment, repulsion, etc. give rise to large scale collective behavior. 
Emergence is physically described as a phase transition to an ordered state (see e.g. the review \cite{Vicsek_Zafeiris_PhysRep12}). Mathematical models often rely on a particle (agent) based description \cite{Aoki1982, Chate2008, Couzin2002, Reynold1987}. Among these models, the time dependent Vicsek model \cite{Vicsek1995} and variants of it have received much attention owing to its structural simplicity, while convincingly describing a wide range of phenomena. The general idea of the Vicsek model is to describe the agents as particles, moving with constant speed. At each time step each particle adjusts its direction towards the mean direction of the particles in its neighborhood subject to some random perturbation.\\

Whilst particle models are a good tool to reproduce and study the observed phenomena, the high number of agents involved make simulations computationally costly considering the amount of detail one is interested in. Hydrodynamic (macroscopic) models are typically easier to simulate and analyze and thereby offer a powerful alternative. Therefore a lot of effort has been made to formulate macroscopic equations starting from particle models \cite{Baskaran2008, Bertin2009, Carrillo2009, Chuang2007, Degond2013, Degond2008, Frouvelle2012, Mogilner1999}. Mathematically one of the main challenges is the lack of conserved quantities, which is typically a key ingredient in the computation of the macroscopic limit (see again \cite{Vicsek_Zafeiris_PhysRep12}). This lack of conservation is frequent in biological systems, which, as opposed to many physical systems, e.g. don't conserve momentum. One strategy to overcome this is the use of \emph{generalized collision invariants (GCI)}, introduced for the first time in \cite{Degond2008} to derive a hydrodynamic model of the Vicsek dynamics (later on referred to as the Self-Organized Hydrodynamic (SOH) model). The passage from individual-based to hydrodynamic models involves the construction of an intermediate kinetic model. For classical physical systems,  momentum or energy conservation translates at the kinetic level into the concept of collision invariants. This is generalized to the GCI concept if the underlying system does not satisfy these conditions. The rigorous passage from kinetic to SOH models has been achieved in \cite{Jiang_Xiong_Zhang_2015}. The GCI concept has been used in a variety of cases \cite{Degond_etal_2015, Degond_Dimarco_Mac_2014, Degond_Liu_2012, Frouvelle2012} and even applied to a model of economics \cite{Degond_Liu_Ringhofer_2014}. Modal analysis of the SOH model can be found in \cite{Degond_Hui_2015}. One important feature of the Vicsek model is the emergence of phase transitions. These are studied in the kinetic of macroscopic framework in \cite{Barbaro_Degond_2014, Bertin2009, Degond_Frouvelle_Liu_2013, Degond_Frouvelle_Liu_2015, Toner_Tu_1995}. Here, we concentrate on the derivation of macroscopic models and postpone the analysis of phase transitions to future works. \\

In this paper we want to apply this method to nematic alignment of polar, rod-shaped particles moving with constant speed in two space dimensions. Nematic alignment describes the situation where upon meeting, two particles either align if their mutual angle is acute or anti-align if their angle is obtuse. This nematic alignment rule typically results from volume exclusion and can be understood as an inelastic collision between polar particles. It differs from polar alignment of polar particles (which is akin to ferromagnetic interactions in spin systems, as is the case of the Vicsek model) and nematic alignment of apolar particles (active nematics, which is typically observed in liquid crystals \cite{DeGennes_Prost_1993},  in fibrous tissues \cite{Degond_Delebecque_Peurichard_2015} or in systems of interacting disks \cite{Degond_Navoret_2015}). Large collections of polar particles which nematically align have been observed to exhibit interesting macroscopic patterns such as high density traveling bands and their properties differ from the other types of alignment \cite{Ginelli2010, Ngo2012}. Therefore the need for a hydrodynamic model derived from particle dynamics arises and this work suggests how to derive such a model. In the derivation we account for purely local interactions as well as weakly non-local interactions which leads to an additional diffusion term in the macroscopic limit.\\

To give an example how the nematic SOH model can be applied to biological systems, we consider the case of myxobacteria. Due to their ability to produce several macroscopic structures, such as waves, spirals and clusters \cite{Dworkin1996, Kaiser2003, Welch2001}, they have already received a lot of attention in terms of mathematical models \cite{Borner2002, Igoshin2001, Igoshin2004, Wu2009}. Since this type of bacteria can internally reverse its polarity, the model has to be complemented with a reaction term, describing these reversals. We demonstrate how such a term can be incorporated in a systematic way into the original model.\\

The rest of the document is structured as follows: In Section \ref{sec:particle_mf}  we describe the particle model for nematic alignment and the resulting mean-field model. The derived alignment operator $Q^0_\text{al}$ is examined in Section \ref{sec:prop_Qal}. Whilst it shares some properties with its polar alignment counterpart, there are also some differences, e.g. the set of equilibria has dimension three (instead of two). Mathematically more unexpected, showing the existence and uniqueness of solutions of the related linear operator, which amounts to proving a Poincar\'e inequality, is slightly more complex than in the polar case \cite{Degond2008}. The section finishes with explicitly characterizing the equilibria and generalized collision invariants of $Q^0_\text{al}$. With these ingredients we derive the macroscopic model in Section \ref{sec:hydro} and prove its hyperbolicity. In Section \ref{sec:nonlocal} we assume weakly non-local interactions and derive the corresponding macroscopic model. Finally in Section \ref{sec:myxo} we apply both macroscopic models to the case of myxobacteria by additionally describing density dependent reversal at the particle level. In the macroscopic model this leads to a reaction term which structurally hints at interesting dynamics expected to take place.

\setcounter{equation}{0}
\section{Particle and Mean Field Model}
\label{sec:particle_mf}

\subsection{Particle Model}
\label{subsec:particle}

As a first step towards a macroscopic model, we start with the description of a time-continuous particle (or agent-based) model. The model's structure is similar to the time-continuous version of the Vicsek model presented in \cite{Degond2008} and we will discuss similarities and differences. We image the particles to resemble stiff rods of equal length and describe the $i$-th particle ($i\in \mathcal{N}:=\{1,\ldots,N\}$) at time $t\geq0$ by its center of mass $X_i(t)\in \mathbb{R}^2$ and its orientation $v(\Theta_i(t))$, where we define
$$ v(\theta) = \left( \begin{array}{c} \cos \theta \\ \sin \theta \end{array} \right), \qquad v^\bot(\theta) = \left( \begin{array}{c} - \sin \theta \\ \cos \theta \end{array} \right)$$ 
as the polar vectors associated with the angle~$\theta$ in a fixed reference frame and where $\Theta_i(t)$ is defined modulo $ 2 \pi$. All particles are assumed to move in the direction of their orientation with constant speed, $v_0$. In particular this means that they cannot drift side-ways (as opposed to the model in \cite{Baskaran2010}). 
Under these hypotheses, the motion of the particles is supposed to  satisfy the following system of stochastic differential equations
\begin{eqnarray}
&&\hspace{-1cm} 
\frac{\dd X_i}{\dd t}=v_0 \, v(\Theta_i(t)), \label{eqn:particle_1} \\
&&\hspace{-1cm} 
\dd \Theta_i= - \nu \, \mbox{Sign} \big( \cos (\Theta_i - \bar \Theta_i) \big) \, \sin (\Theta_i - \bar \Theta_i) \,  \dd t 
+ \big( 2\,d \, \cos^2 (\Theta_i - \bar \Theta_i) \big)^{1/2} \dd B^i_t .  \label{eqn:particle_2}
\end{eqnarray}
Here $\bar \Theta_i = \bar \Theta_i(t)$ is the average nematic alignment direction. It is defined as follows: First define a mean 'nematic' current
$$ J_i(t)=\sum_{k\in\mathcal{N}: |X_k(t)-X_i(t)|\leq R} v\big(2\Theta_k(t)\big), $$ 
where $R>0$ is the interaction radius. Next define $\bar \Theta_i(t)$ as an angle of lines (i.e. an angle defined modulo $\pi$)  such that:
\begin{equation}
v\big(2 \bar \Theta_i(t) \big) = \frac{J_i(t)}{|J_i(t)|}. 
\label{eq:def_bar_Theta_i}
\end{equation}
We note that the definition of the vector on the left-hand side of (\ref{eq:def_bar_Theta_i}) does not depend on the choice of the reference frame used to measure the angles $\Theta_i$ as $\bar \Theta_i(t)$ can be equivalently defined by 
$$ \sum_{k\in\mathcal{N}: |X_k(t)-X_i(t)|\leq R}  \sin 2 \big(\Theta_k(t) - \bar \Theta_i(t) \big)   = 0, $$ 
and we notice that this relation is invariant under a uniform translation of all angles. To obtain a unique representation of the angles, we will use the convention that the particle angles $\Theta_i(t)$ are taken in $[-\pi, \pi)$ and the nematic angles $\bar \Theta_i(t)$ are taken in $[-\pi/2,\pi/2)$. For simplicity, we will identify $[-\pi, \pi)$ with ${\mathbb R}/2 \pi {\mathbb Z}$ and with the unit circle ${\mathbb S}^1$. The notation $\!\dd B_t^i$ describes independent Brownian motions with intensity $d\,\cos^2 (\Theta_i - \bar \Theta_i) $ with $d > 0$. We are aiming for a separation of the particles into two (locally defined) groups, one with $\cos (\Theta_i - \bar \Theta_i)>0$ and one with $\cos (\Theta_i - \bar \Theta_i)<0$. Therefore, diffusion must vanish near $\cos (\Theta_i - \bar \Theta_i)=0$, to avoid particles crossing the boundary between these two angular regions. The factor $\cos^2 (\Theta_i - \bar \Theta_i)$ is intended to avoid that the two groups mix. The first term in the velocity evolution describes alignment with the bacteria within the interaction radius  $R> 0$ with an alignment frequency $\nu \, \mbox{Sign} (\cos (\Theta_i - \bar \Theta_i) ) \, \sin (\Theta_i - \bar \Theta_i)$ with $\nu\geq 0$. Due to the factor $\mbox{Sign} (\cos (\Theta_i - \bar \Theta_i) )$, alignment occurs parallel or anti-parallel to the mean angle $\bar \Theta_i$, depending on whether $\Theta_i$ is closer to $\bar \Theta_i$ or $\bar \Theta_i + \pi$. Note that $\mbox{Sign} (\cos (\Theta_i - \bar \Theta_i) ) \, \sin (\Theta_i - \bar \Theta_i)$ is a $\pi$-periodic function and so, does not depend on the choice of $\bar \Theta_i$ or $\bar \Theta_i + \pi$ as the representative of the corresponding angle of lines. 

\begin{remark} 
Also other notions of nematic alignment between several particles are possible. In \cite{Peruani2011} the nematic mean direction was defined as
\begin{align}
\label{eqn:Chate}
\bar\Theta_i(t) = \arg\left(\sum_{k\in\mathcal{N}:|X_k(t)-X_i(t)|\leq R}{\rm Sign}(\cos(\Theta_i(t)-\Theta_k(t)))v(\Theta_k)\right).
\end{align}
The main difference between this definition and \eqref{eq:def_bar_Theta_i} is that the nematic angle defined in  \eqref{eq:def_bar_Theta_i} depends on the particle index $i$ only via the position $X_i(t)$, but is independent of the direction of the $i$-th particle, $\Theta_i$. The outcome of \eqref{eqn:Chate} on the other hand depends on both the position and the direction of the $i$-th particle. Whilst the definition from \cite{Peruani2011} might be closer to reality, if the dispersion around the mean angle is small, both types of nematic alignment lead to the same behavior. However, the independence from $\Theta_i$ makes  \eqref{eq:def_bar_Theta_i} mathematically much easier to handle.
\end{remark}

\subsection{Mean Field Model}
\label{ssec:meanfield}

We follow the usual procedure to derive the mean field limit as described for example in \cite{Degond2008}. The procedure starts with the definition of the empirical distribution function $f^N(x,\theta,t)$, with $x\in \mathbb{R}^2$, $\theta \in {\mathbb R}$ (modulo $2 \pi$) and $t>0$, according to 
$$ f^N(x,\theta,t) = \frac{1}{N} \sum_{i \in {\mathcal N}} \delta_{(X_i(t), \Theta_i(t))} (x, \theta),$$
where $\delta_{(X_i(t), \Theta_i(t))} (x, \theta)$ is the Dirac delta in $(x,\theta)$ space located at $(X_i(t), \Theta_i(t))$ at time~$t$. In the absence of noise ($d=0$), $ f^N$ is a deterministic measure that  satisfies a first order kinetic equation if we assume that all $(X_i, V_i)$ fulfill \eqref{eqn:particle_1}-\eqref{eq:def_bar_Theta_i}. If $d \not = 0$, $f^N$ is a random measure whose law satisfies a second order kinetic equation. 
Then as  $N\rightarrow \infty$, $f^N$ approximates a distribution function $f(x,\theta,t)$, satisfying the following Kolmogorov-Fokker-Planck type equation:
\begin{equation}
\partial_t f+ v_0\nabla_x\cdot (v(\theta) \, f)=\tilde Q_\text{al}f, \label{eqn_meanfield}
\end{equation}
where the collision operator $\tilde Q_\text{al}$ is given by
\begin{eqnarray}
&&\hspace{-1cm}
\tilde Q_\text{al}f=\partial_\theta  \Big[ \nu \, \mbox{Sign} \big(\cos(\theta- \bar \Theta_f) \big) \, \sin (\theta- \bar \Theta_f) \, f 
+ d \, \cos^2(\theta- \bar \Theta_f) \, \partial_\theta f \Big], 
\label{eqn:scaling_Qal} 
\end{eqnarray}
and where $ \bar \Theta_f$ is given by
\begin{eqnarray}
&&\hspace{-1cm}
v \big( 2 \bar \Theta_f(x,t) \big)  = \frac{J_f(x,t)}{|J_f(x,t)|}, 
\label{eqn:scaling_jf_0} \\
&&\hspace{-1cm}
J_f(x,t)=\int_{|y-x|\leq R} \int_{-\pi}^{\pi} v( 2 \theta)\, f(y,\theta,t) \dd y \dd \theta. 
\label{eqn:scaling_jf} 
\end{eqnarray}
We note from (\ref{eqn:scaling_jf_0}) that $\bar \Theta_f$ is an angle of lines, defined modulo $\pi$ and a representative of $\bar \Theta_f$ will be chosen in $[-\pi/2, \pi/2)$. We also note that (\ref{eqn:scaling_jf_0})-(\ref{eqn:scaling_jf}) can be equivalently written
$$\int_{|y-x|\leq R} \int_{-\pi}^{\pi} 
\sin 2 \big( \theta - \bar \Theta_f(x,t) \big) \, f(y,\theta,t) \dd y \dd \theta = 0. $$
At this point the alignment interaction is non-local which is reflected in the fact that the computations for $J_f(x,t)$ involve integrals over a ball centered around~$x$. As a next step we introduce dimensionless variables. Since we are interested in macroscopic phenomena, we use a hydrodynamic scaling. In the following we will do the non-dimensionalization and the hydrodynamic scaling in one step. A more detailed description for a similar equation can be found in \cite{Degond2013}. On a microscopic scale we use the reference time $t_0=1/\nu$ and reference length $x_0=v_0\, t_0$. We also introduce a scaled diffusion constant $D=d\,t_0$. On the macroscopic scale we want to use units $t_0'=\frac{t_0}{\e}$ and $x_0'=\frac{x_0}{\e}$, which are large compared to the microscopic ones. We introduce the new (dimensionless, macroscopic) variables $\hat x=\frac{x}{x_0'}=\e \frac{x}{x_0}$ and $\hat t=\frac{t}{t_0'}=\e \frac{t}{t_0}$ and set $\hat f(\hat x, \theta, \hat t)=\frac{f(x,\theta,t)}{(1/x_0')^2}=\left(\frac{x_0}{\e}\right)^2 f(x,\theta,t)$. Additionally we scale $\hat R=\frac{R}{x_0'}=\e \frac{R}{x_0} = \varepsilon r$ with $r= \frac{R}{x_0}  = {\mathcal O}(1)$. Note that the scaling of $\hat R$ supposes that the interactions between the particles are local. Another scaling which takes better account of the non-locality of the interactions will be investigated later on (see Section \ref{sec:nonlocal}). With this scaling and dropping the hats for better readability we find that $f=f^\varepsilon$ satisfies the problem:
\begin{eqnarray}
&&\hspace{-1cm}
\varepsilon \big( \partial_t f^\varepsilon+ \nabla_x\cdot (v(\theta) \, f^\varepsilon) \big)=\tilde Q^\varepsilon_\text{al} f^\varepsilon, \label{eqn_meanfield_scaled} \
\end{eqnarray}
with (omitting the variables $x$ and $t$):
\begin{eqnarray}
&&\hspace{-1cm}
\tilde Q^\varepsilon_\text{al}f(\theta)=\partial_\theta  \Big[ \mbox{Sign} \big(\cos(\theta- \bar \Theta_f^\varepsilon) \big) \, \sin (\theta- \bar \Theta_f^\varepsilon) \, f (\theta)
+ D \, \cos^2(\theta- \bar \Theta_f^\varepsilon) \, \partial_\theta f (\theta)\Big], 
\label{eqn:Qal_scaled} 
\end{eqnarray}
and 
\begin{eqnarray}
&&\hspace{-1cm}
v \big( 2 \bar \Theta^\varepsilon_f(x,t) \big)  = \frac{J^\varepsilon_f (x,t)}{|J^\varepsilon_f (x,t)|}, 
\label{eqn:Thetaf_scaled} \\
&&\hspace{-1cm}
J_f^\varepsilon(x,t)= \int_{-\pi}^{\pi} \int_{|y-x|\leq \varepsilon r} v(2 \theta)  \, f(y,\theta,t) \dd y \dd \theta. 
\label{eqn:scaling_jf_2} 
\end{eqnarray}
By Taylor expansion we get:
$$
\bar \Theta^\varepsilon_f(x,t)=\bar \theta_f(x,t)+\mathcal{O}(\e^2),
$$
where
\begin{eqnarray}
&&\hspace{-0.5cm}
v \big( 2 \bar \theta_f(x,t) \big) = \frac{j_f(x,t)}{|j_f(x,t)|}, \label{eqn:Omega_0} \\
&&\hspace{-0.5cm}
j_f(x,t)=\int_{-\pi}^{\pi} v(2\theta) \, f(x,\theta,t) \dd \theta. 
\label{eqn:Omega} 
\end{eqnarray}
Again, $\bar \theta_f$ is an angle of lines, defined modulo $\pi$ and a representative of $\bar \theta_f$ is chosen in $[-\pi/2, \pi/2)$. As above,  (\ref{eqn:Omega_0})-(\ref{eqn:Omega}) can be equivalently written
$$\int_{-\pi}^{\pi} \sin 2 \big( \theta - \bar \theta_f(x,t) \big) \, f(x,\theta,t) \dd \theta = 0. $$

Dropping the terms of order $\varepsilon^2$ or smaller, we  find that $f = f^\varepsilon(x,\theta,t)$ satisfies the perturbation problem: 
\begin{eqnarray}
&&\hspace{-1cm}
\e\big( \partial_t f^\varepsilon+\nabla_x \cdot (v(\theta) \,f^\varepsilon) \big)=Q^0_\text{al} f^\varepsilon, \label{eqn:scaling_hydro} \\
&&\hspace{-1cm}
Q^0_\text{al}f=\partial_\theta  \Big[ \mbox{Sign} \big(\cos(\theta- \bar \theta_f) \big) \, \sin (\theta- \bar \theta_f) \, f 
+ D \, \cos^2(\theta- \bar \theta_f) \, \partial_\theta f \Big], 
\label{eqn:Qal}
\end{eqnarray}
where $\bar \theta_f$ is given by \eqref{eqn:Omega_0}.



\setcounter{equation}{0}
\section{Properties of the alignment operator $Q^0_\text{al}$}
\label{sec:prop_Qal}

\subsection{Functional setting}
\label{subsec_funct}

Let $\kappa = 1/D$. We can write:
\begin{eqnarray}
&&\hspace{-1cm}
Q^0_\text{al}f= D \partial_\theta  \Big\{ \cos^2(\theta- \bar \theta_f) \Big[ \kappa \frac{\mbox{Sign} \big(\cos(\theta- \bar \theta_f) \big) \, \sin (\theta- \bar \theta_f)}{\cos^2(\theta- \bar \theta_f)} \, f + \partial_\theta f \Big] \Big\}. 
\label{eqn:Qal_shape}
\end{eqnarray}
Now, we introduce 

\begin{definition}
For any angle of lines $\theta_0 \in [- \pi/2,\pi/2)$, we define the Generalized von Mises (GVM) distribution with direction $\theta_0$ and concentration parameter $\kappa$ as:

\begin{equation}
M_{\theta_0} (\theta) = \frac{1}{Z_\kappa} \exp \Big( - \frac{\kappa}{|\cos (\theta - \theta_0)|} \Big) , \quad \theta \in  [- \pi,\pi), \label{eqn:M_Omega}
\end{equation}
where 
$$ Z_\kappa = \int_{\cos \theta >0}  \exp \Big( - \frac{\kappa}{\cos \theta}
\Big) \dd \theta, $$
is the normalization coefficient such that 
\begin{equation} 
\int_{\cos \theta >0} M_{\theta_0} (\theta) \dd \theta = \int_{\cos \theta <0} M_{\theta_0} (\theta) \dd \theta = 1.
\label{eq:GVM_norm}
\end{equation}
\label{def:GVM}
\end{definition}

\begin{figure}
\center
\includegraphics[width=0.5\textwidth]{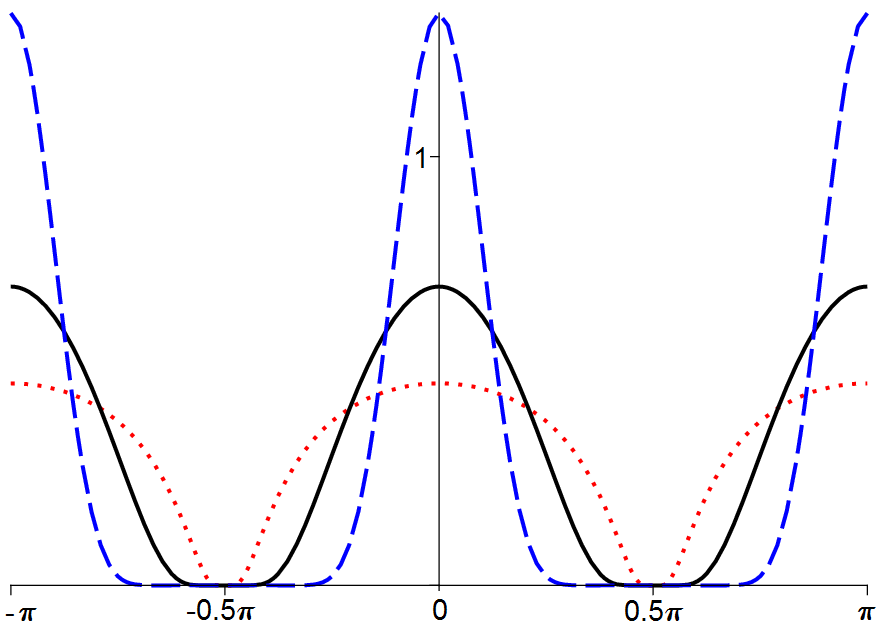}
\caption{$M_0(\theta)$ for $\kappa=0.5, 2,10$ (red-dotted, black-solid, blue-dashed).}
\label{fig:GVM}
\end{figure}

\medskip
\noindent
Figure \ref{fig:GVM} shows $M_0(\theta)$ for different values of $\kappa$. With this definition, we have

\begin{lemma}
\label{lem:Qshape}
The alignment operator $Q^0_\text{al}$ can be written:
\begin{equation}
Q^0_\text{al}f (\theta) = D \partial_\theta  \Big[ \cos^2(\theta- \bar \theta_f) \, M_{\bar \theta_f}  \, \partial_\theta \Big( \frac{f}{M_{\bar \theta_f}} \Big)  \Big]. 
\label{eq:Qal_3}
\end{equation}
\end{lemma}

\medskip 
\noindent 
{\bf Proof.}
Using \eqref{eqn:M_Omega} we calculate
\begin{eqnarray*}
M_{\bar \theta_f} \, \partial_\theta \Big( \frac{f}{M_{\bar \theta_f}} \Big) &=&\partial_\theta f - f \, \partial_\theta (\log M_{\bar \theta_f} ) \\
&=&\partial_\theta f + \kappa f \, \partial_\theta \left(\frac{1}{|\cos (\theta - \bar \theta_f) |}\right)\\
&=& \partial_\theta f + \kappa f \, \frac{\mbox{Sign} (\cos (\theta - \bar \theta_f)) \, \sin (\theta - \bar \theta_f)}{\cos^2 (\theta - \bar \theta_f)}, 
\end{eqnarray*}
which, upon multiplication by $\cos^2(\theta- \bar \theta_f)$, yields \eqref{eqn:Qal_shape} and shows the claim. \endproof

\medskip
\noindent
We now give a more precise functional setting in which to define the collision operator $Q_\text{al}^0$. For this purpose, we introduce the following:  

\begin{definition}
(i) Let $\theta_0 \in [-\pi/2, \pi/2)$ be a given angle of lines. The linear operator ${\mathcal L}_{\theta_0}$ is defined as
\begin{equation} 
{\mathcal L}_{\theta_0} \varphi = D \, \frac{1}{M_{\theta_0}} \, \partial_\theta  \big( \cos^2(\theta- \theta_0) \, M_{\theta_0}  \, \partial_\theta \varphi  \big), 
\label{eqn:L}
\end{equation}
for any smooth function $\varphi$ defined on ${\mathbb R}/(2 \pi {\mathbb Z})$ with values in ${\mathbb R}$. \\
(ii) The space ${\mathcal H}_{\theta_0}$ is defined by 
$$ {\mathcal H}_{\theta_0} = \Big\{ g: \, {\mathbb R}/2\pi {\mathbb Z} \to {\mathbb R}, \mbox{ measurable, such that } \int_{-\pi}^{\pi} |g(\theta)|^2 \, M_{\theta_0} (\theta) \dd \theta < \infty \Big\}. $$
It is a Hilbert space when endowed with the inner product 
$$(g,h)_{{\mathcal H}_{\theta_0}} = \int_{-\pi}^{\pi} g(\theta) \, h(\theta) \, M_{\theta_0} (\theta) \dd \theta,  $$
and associated norm $ \|g\|_{{\mathcal H}_{\theta_0}}^2 = (g,g)_{{\mathcal H}_{\theta_0}}$. \\
(iii) The space ${\mathcal V}_{\theta_0}$ is defined by 
$$ {\mathcal V}_{\theta_0} = \Big\{ g \in  {\mathcal H}_{\theta_0}, \, \, \mbox{ such that } \,  \, \sqrt{M_{\theta_0} (\theta)} \, |\cos (\theta - \theta_0)| \, \partial_\theta g(\theta) \in L^2(-\pi, \pi)  \Big\}, $$
where $\partial_\theta g$ is defined in the distributional sense. It is a Hilbert space when endowed with the inner product 
\begin{eqnarray*}
&&\hspace{-1cm}
(g,h)_{{\mathcal V}_{\theta_0}} = (g,h)_{{\mathcal H}_{\theta_0}} \\
&&\hspace{1cm}
+ \int_{-\pi}^{\pi} \Big( \sqrt{M_{\theta_0} (\theta)} \, |\cos (\theta - \theta_0)| \, \partial_\theta g(\theta)  \Big) \Big( \sqrt{M_{\theta_0} (\theta)} \, |\cos (\theta - \theta_0)| \, \partial_\theta h(\theta)  \Big) \dd \theta, 
\end{eqnarray*}
and associated norm $\|g\|_{{\mathcal V}_{\theta_0}}^2 = (g,g)_{{\mathcal V}_{\theta_0}}$. \\
(iv) ${\mathcal L}_{\theta_0}$ defines a continuous mapping from ${\mathcal V}_{\theta_0}$ into the dual space $ {\mathcal V}_{\theta_0}'$ of  ${\mathcal V}_{\theta_0}$ by: 
\begin{eqnarray}
&&\hspace{-1.5cm} 
\langle {\mathcal L}_{\theta_0} \varphi, \psi \rangle_{{\mathcal V}_{\theta_0}',{\mathcal V}_{\theta_0}}  \nonumber \\
&&\hspace{-1cm} 
= - D \int_{-\pi}^{\pi} \big( \sqrt{M_{\theta_0} (\theta)} \, |\cos (\theta - \theta_0)| \, \partial_\theta \varphi \big) \, \big( \sqrt{M_{\theta_0} (\theta)} \, |\cos (\theta - \theta_0)| \, \partial_\theta  \psi \big)\dd \theta, 
\label{eq:Lweak}
\end{eqnarray}
for any $\varphi$, $\psi \in {\mathcal V}_{\theta_0}$. Here,  
$\langle \cdot , \cdot \rangle_{{\mathcal V}_{\theta_0}',{\mathcal V}_{\theta_0}}$ denotes the duality bracket between ${\mathcal V}_{\theta_0}$ and its dual. This duality bracket extends the inner product of ${\mathcal H}_{\theta_0}$, i.e.\ for $f \in  {\mathcal H}_{\theta_0} \subset {\mathcal V}_{\theta_0}'$ and $g \in {\mathcal V}_{\theta_0}$, we have  $\langle f , g \rangle_{{\mathcal V}_{\theta_0}',{\mathcal V}_{\theta_0}} = \int f \, g \, M_{\theta_0} \dd \theta$. The operator ${\mathcal L}_{\theta_0} $ is formally self adjoint, i.e.\ 
$\langle {\mathcal L}_{\theta_0} \varphi, \psi \rangle_{{\mathcal V}_{\theta_0}',{\mathcal V}_{\theta_0}} = \langle {\mathcal L}_{\theta_0} \psi, \varphi \rangle_{{\mathcal V}_{\theta_0}',{\mathcal V}_{\theta_0}}$, for all $\varphi$, $\psi \in {\mathcal V}_{\theta_0}$.
\\
(v) The operator $Q^0_\text{al}$ is defined for functions $f \in L^1(-\pi,\pi)$ (which guarantees that $\bar \theta_f$ is defined) such that $f/M_{\bar \theta_f} \in {\mathcal V}_{\bar \theta_f}$ by:
\begin{equation}
Q^0_\text{al} f =  M_{\bar \theta_f} {\mathcal L}_{\bar \theta_f} \Big( \frac{f}{M_{\bar \theta_f}} \Big). 
\label{eq:defQal}
\end{equation}
\label{def:linear_op}
\end{definition}

\noindent
By abuse of notation, we will simply write (\ref{eq:Lweak}) as follows: 
\begin{eqnarray*}
&&\hspace{-1.5cm} 
\langle {\mathcal L}_{\theta_0} \varphi , \psi \rangle_{{\mathcal V}_{\theta_0}',{\mathcal V}_{\theta_0}} 
= - D \int_{-\pi}^{\pi} \partial_\theta \varphi  \, \,  \partial_\theta  \psi \, M_{\theta_0} (\theta) \, \cos^2 (\theta - \theta_0) \dd \theta. 
\end{eqnarray*}
Now we investigate the solvability of the equation 
\begin{equation} 
{\mathcal L}_{\theta_0} \varphi = \psi, 
\label{eq:Lf=h}
\end{equation}
for a given function $\psi$. In view of (\ref{eq:Lweak}), this problem can be recast into the following variational formulation: 
\begin{eqnarray}
&& \hspace{-1cm}
\varphi \in {\mathcal V}_{\theta_0}, \label{eq:vf1} \\
&& \hspace{-1cm}
- D \int_{-\pi}^{\pi} \partial_\theta \varphi \, \, \partial_\theta \chi \, M_{\theta_0} \, \cos^2(\theta - \theta_0) \dd \theta = \int_{-\pi}^{\pi} \psi \, \chi \, M_{\theta_0} \dd \theta, \qquad \forall \chi \in {\mathcal V}_{\theta_0}. \label{eq:vf2}
\end{eqnarray}
For this variational formulation, we have the following

\begin{proposition}
(i) (Homogeneous problem) Suppose $\psi = 0$. Then, the unique solutions of (\ref{eq:vf1}), (\ref{eq:vf2}) are piecewise constant functions, i.e. there exist $C_\pm \in {\mathbb R}$ such that 
\begin{equation}
\varphi(\theta) = \left\{ \begin{array}{lll} C_+ & \mbox{ if } & \cos (\theta - \theta_0) >0, \\
 C_- & \mbox{ if } & \cos (\theta - \theta_0) <0. \\
\end{array} \right. 
\label{eq:equil}
\end{equation}
(ii) (Nonhomogeneous problem) Let $\psi \in {\mathcal V}_{\theta_0}$. Then, there exists a solution to
(\ref{eq:vf1}), (\ref{eq:vf2}) if and only if the following two solvability conditions hold simultaneously: 
\begin{equation}
\int_{\cos (\theta - \theta_0) >0} \psi \, M_{\theta_0} \dd \theta = 0, \quad  \int_{\cos (\theta - \theta_0) <0} \psi \, M_{\theta_0} \dd \theta = 0.
\label{eq:solv}
\end{equation}
Under these solvability conditions, there exists a unique solution $\varphi$ to (\ref{eq:vf1}), (\ref{eq:vf2}) which satisfies the additional conditions 
\begin{equation}
\int_{\cos (\theta - \theta_0) >0} \varphi \, M_{\theta_0} \dd \theta = 0, \quad  \int_{\cos (\theta - \theta_0) <0} \varphi \, M_{\theta_0} \dd \theta = 0,
\label{eq:cancel}
\end{equation}
and the general solution in ${\mathcal V}_{\theta_0}$ to (\ref{eq:vf1}), (\ref{eq:vf2}) is the sum of the unique solution satisfying (\ref{eq:cancel}) and the general solution (\ref{eq:equil}) of the homogeneous problem. \\
(iii) (Invariance) Let $\theta_0 = 0$, $\psi_0 \in {\mathcal V}_0$ satisfying the solvability condition (\ref{eq:solv}) and $\varphi_0$ be the unique solution to (\ref{eq:vf1}), (\ref{eq:vf2}) satisfying (\ref{eq:cancel}). Let now $\theta_0 \in [-\pi/2,\pi/2)$, $\theta_0 \not = 0$ and define $\psi$ and $\varphi$ such that $\psi (\theta) = \psi_0(\theta - \theta_0)$, $\varphi(\theta) = \varphi_0(\theta-\theta_0)$. Then, $\varphi$ is the unique solution to (\ref{eq:vf1}), (\ref{eq:vf2}) satisfying (\ref{eq:cancel}) with that value of $\theta_0$.
\label{prop:lineareq}
\end{proposition}

\medskip
\noindent
{\bf Proof.} We easily verify that we can choose $\theta_0 = 0$ and thus we drop the indices $\theta_0$ in the definitions of the operators and spaces. 

\medskip
\noindent
(i) We first verify that the function $\varphi$ given by (\ref{eq:equil}) belongs to ${\mathcal V}$ and satisfies (\ref{eq:vf1}), (\ref{eq:vf2}) with $\psi = 0$. We clearly have $\varphi \sqrt{M} \in L^2(-\pi,\pi)$ and
$$\partial_\theta \varphi = (C_+-C_-) \, ( \delta_{-\pi/2} - \delta_{\pi/2}), $$ 
where $\delta_{\pm \pi/2}$ is the Dirac delta at $\theta = \pm \pi/2$. Since $|\cos \theta|\, \sqrt{M}$ is a $C^\infty$ function which takes the value $0$ at $\theta = \pm \pi/2$, we have $|\cos \theta| \, \sqrt{M} \, \partial_\theta \varphi = 0 \in L^2(-\pi,\pi)$ and $\varphi$ satisfies~(\ref{eq:vf2}) with right-hand side equal to zero.  Now, we suppose that $\varphi \in {\mathcal V}$ is a function satisfying~(\ref{eq:vf2}) with right-hand side equal to zero. Taking $\chi = \varphi$ as a test function in~(\ref{eq:vf2}), we get 
$|\cos \theta| \, \sqrt{M} \, \partial_\theta \varphi = 0$. This implies that $\varphi$ is constant on all connected components of the set $(|\cos \theta| \, \sqrt{M})^{-1}({\mathbb R})$ and is $2 \pi$-periodic, which implies that $\varphi$ is of the form (\ref{eq:equil}). 

\medskip
\noindent
(ii) Let $\psi \in {\mathcal V}$ be such that there exists a solution $\varphi$ of the variational formulation  (\ref{eq:vf1}),  (\ref{eq:vf2}). Taking $\chi$ of the form (\ref{eq:equil}), we get that $\psi$ must satisfy the solvability conditions~(\ref{eq:solv}). Now, we denote by $\tilde {\mathcal V}$ the subspace of ${\mathcal V}$ defined by 
$$ \tilde {\mathcal V} = \Big\{ \varphi \in {\mathcal V} \mbox{ such that } \int_{\cos \theta >0} \varphi \, M \dd \theta = 0 \mbox{ and } \int_{\cos \theta <0} \varphi \, M \dd \theta = 0 \Big\}. $$
Suppose that $\psi \in {\mathcal V}$ satisfies (\ref{eq:solv}) and that $\varphi$ is a solution of the following variational formulation: 
\begin{eqnarray}
&& \hspace{-1cm}
\varphi \in \tilde {\mathcal V}, \label{eq:tvf1} \\
&& \hspace{-1cm}
D \int_{-\pi}^{\pi} \partial_\theta \varphi \, \, \partial_\theta \chi \, M \, \cos^2\theta \dd \theta = - \int_{-\pi}^{\pi} \psi \, \chi \, M \dd \theta, \qquad \forall \chi \in \tilde {\mathcal V}. \label{eq:tvf2}
\end{eqnarray}
Then, we claim that $\varphi$ is a solution of (\ref{eq:vf1}),  (\ref{eq:vf2}). Indeed, let $\chi \in {\mathcal V}$. Define $\tilde \chi = \chi - \chi_0$ with $\chi_0$ constructed according to (\ref{eq:equil}) with the constants:
$$ C_\pm = \int_{\pm \cos \theta >0} \chi \, M \dd \theta. $$
Then, $\tilde \chi \in \tilde {\mathcal V}$. Therefore, we can apply (\ref{eq:tvf2}) with test function $\tilde \chi$. But since $\chi_0$ is of the form (\ref{eq:equil}) and since $\psi$ satisfies (\ref{eq:solv}), the contribution of $\chi_0$ to both the left and right hand sides of (\ref{eq:tvf2}) vanish. Consequently, $\varphi$ satisfies (\ref{eq:vf1}),  (\ref{eq:vf2}) as claimed. 

\medskip
\noindent
We now show that the variational formulation (\ref{eq:tvf1}), (\ref{eq:tvf2}) admits a unique solution. For this purpose, we apply Lax-Milgram's Lemma. First, we note that $\tilde {\mathcal V}$ is a closed subspace of ${\mathcal V}$. Indeed, it is obvious that the linear forms $\varphi \to \int_{\pm \cos \theta >0} \varphi \, M \dd \theta$ are continuous. Therefore, $\tilde {\mathcal V}$ is a Hilbert space. Then, it is clear that the left-hand side of (\ref{eq:tvf2}) defines a symmetric continuous bilinear form $a(\varphi,\chi)$ on $\tilde {\mathcal V}$ and that the right-hand side of (\ref{eq:tvf2}) defines a continuous linear form $L(\chi)$ on $\tilde {\mathcal V}$. To apply Lax-Milgram's Lemma, it remains to prove that $a$ is coercive. This amounts to proving a Poincar\'e inequality, i.e.\ that there exists a constant $C>0$ such that for all $\varphi \in \tilde {\mathcal V}$ we have:
\begin{equation}
\int_{-\pi}^{\pi} \big| \partial_\theta \varphi\big|^2  \, M \, \cos^2 \theta \dd \theta \, \, \geq \, \, C \,  \int_{-\pi}^{\pi} | \varphi |^2  \, \, M \dd \theta. 
\label{eq:poincare}
\end{equation}
We show that 
\begin{equation}
\int_{\cos \theta >0} \big| \partial_\theta \varphi\big|^2  \, M \, \cos^2 \theta \dd \theta \, \, \geq \, \, C \,  \int_{\cos \theta >0} | \varphi |^2  \, \, M \dd \theta. 
\label{eq:poincare_2}
\end{equation}
Obviously, if such a result is true, it is also true if the integration domains are changed to the set $\{ \cos \theta <0 \}$ in both integrals, and this implies (\ref{eq:poincare}). We now concentrate on obtaining (\ref{eq:poincare_2}). Let $\psi \in \tilde {\mathcal V}$ and let $\theta_1$ and $\theta_2$ be such that $\cos \theta_1 >0$, $\cos \theta_2 >0$. We have:
$$ \psi(\theta_2) - \psi(\theta_1) = \int_{\theta_1}^{\theta_2} \partial_\theta \psi (\theta) \dd \theta. $$
Multiplying by $M(\theta_1)$, integrating with respect to $\theta_1$ on $(-\pi/2,\pi/2)$ and using that 
$\int_{\cos \theta_1 >0} \psi(\theta_1) \, M(\theta_1) \dd \theta_1 = 0$ and $\int_{\cos \theta_1 >0} M(\theta_1) \dd \theta_1 = 1$, we get:
$$ \psi(\theta_2) = \int_{\cos \theta_1 >0}  \Big( \int_{\theta_1}^{\theta_2} \partial_\theta \psi (\theta) \dd \theta \Big) \, M(\theta_1) \dd \theta_1. $$
Taking the square, multiplying by $M(\theta_2)$ and integrating with respect to $\theta_2$ on $(-\pi/2,\pi/2)$ leads to:
$$ \int_{\cos \theta_2 >0} |\psi(\theta_2)|^2 \, M(\theta_2) \dd \theta_2 = \int_{\cos \theta_2 >0}  \Big| \int_{\cos \theta_1 >0} \Big( \int_{\theta_1}^{\theta_2} \partial_\theta \psi (\theta) \dd \theta \Big)  \, M(\theta_1) \dd \theta_1 \Big|^2  \, M(\theta_2) \dd \theta_2, $$
and by the Cauchy-Schwarz inequality with respect to the inner product on ${\mathcal H}$, again using that $\int_{\cos \theta_1 >0} M(\theta_1) \dd \theta_1 = 1$, we obtain
\begin{eqnarray}
\int_{\cos \theta_2 >0} |\psi(\theta_2)|^2 \, M(\theta_2) \dd \theta_2 &\leq&  \int_{\cos \theta_1 >0} \int_{\cos \theta_2 >0}  \Big| \int_{\theta_1}^{\theta_2} \partial_\theta \psi (\theta) \dd \theta \Big|^2  \, M(\theta_1) \, M(\theta_2) \dd \theta_1 \dd  \theta_2 \nonumber \\
&\leq& I_1 + \ldots + I_4, 
\label{eq:poincare3}
\end{eqnarray}
where $I_1, \ldots I_4$ are the integrals with respect to the sets $\Omega_1$ to $\Omega_4$ with  $\Omega_1 = \{ |\sin \theta_1| < 1/2, \, \, |\sin \theta_2| < 1/2\}$, $\Omega_ 2= \{ |\sin \theta_1| < 1/2, \, \, |\sin \theta_2| > 1/2\}$, $\Omega_3 = \{ |\sin \theta_1| > 1/2, \, \, |\sin \theta_2| < 1/2\}$, $\Omega_4 = \{ |\sin \theta_1| > 1/2, \, \, |\sin \theta_2| > 1/2\}$ respectively. For $|\sin \theta| < 1/2$, there exist constants $C_1$ and $C_2$ such that $0 < C_1 \leq M(\theta) \leq C_2 <\infty$ and $0 < C_1 \leq \cos^2 \theta  \, M(\theta) \leq C_2 <\infty$. So, for $(\theta_1,\theta_2) \in \Omega_1$, using Cauchy-Schwarz inequality, we have: 
\begin{equation}
 \int_{\theta_1}^{\theta_2} \partial_\theta \psi (\theta) \dd  \theta \leq C \Big( \int_{\cos \theta >0} | \partial_\theta \psi (\theta)|^2  \, M (\theta) \, \cos^2 \theta \dd  \theta \Big)^{1/2}, 
\label{eq:inttheta}
\end{equation}
and consequently, 
\begin{equation}
 I_1 \leq C \int_{\cos \theta >0} | \partial_\theta \psi (\theta) |^2  \, M (\theta) \, \cos^2 \theta \dd  \theta , 
\label{eq:intparthetapsi}
\end{equation}
where here and in the following, we denote by $C$ generic constants. 
Take now $(\theta_1,\theta_2) \in \Omega_4$. Then, using Cauchy-Schwarz inequality, we can write: 
\begin{eqnarray}
&&\hspace{-1cm} 
\int_{\theta_1}^{\theta_2} \partial_\theta \psi (\theta) \dd  \theta = 
\int_{\theta_1}^{\theta_2} |\cos \theta| \, \sqrt{M(\theta)} \, \partial_\theta \psi (\theta) \, \big(  |\cos \theta| \, \sqrt{M(\theta)} \big)^{-1} \dd  \theta \nonumber \\
&&\hspace{0cm} 
\leq  \left( \int_{\theta_1}^{\theta_2} | \partial_\theta \psi (\theta) |^2 \, \cos^2 \theta \, M(\theta) \dd  \theta \right)^{1/2} \left( \int_{\theta_1}^{\theta_2}    \frac{d \theta}{\cos^{2} \theta \, M(\theta)}   \right)^{1/2} . 
\label{eq:inttheta2}
\end{eqnarray}
Now, we remark that 
$$ \partial_\theta \left( \frac{1}{M(\theta)} \right) = \kappa \, \frac{1}{M(\theta)} \frac{\sin \theta}{\cos^2 \theta}, $$
and deduce that
\begin{eqnarray}
&&\hspace{-2.1cm}
\int_{\theta_1}^{\theta_2}    \frac{d \theta}{\cos^{2} \theta \, M(\theta)} =
\frac{1}{\kappa} \int_{\theta_1}^{\theta_2}    \partial_\theta \left( \frac{1}{M(\theta)} \right) \, \frac{1}{\sin \theta} \dd  \theta \nonumber\\
&&\hspace{1cm}
= \frac{1}{\kappa} \left( \frac{1}{M(\theta_2) \, \sin \theta_2} - \frac{1}{M(\theta_1) \, \sin \theta_1} \right) + \frac{1}{\kappa}  \int_{\theta_1}^{\theta_2}    \frac{1}{M(\theta)} \, \frac{\cos \theta} {\sin^2 \theta} \dd  \theta.   
\label{eq:inttheta3}
\end{eqnarray}
Then, we have
\begin{eqnarray}
&&\hspace{-1cm}
M(\theta_1) \, M(\theta_2) \, \int_{\theta_1}^{\theta_2}    \frac{\dd \theta}{\cos^{2} \theta \, M(\theta)} \nonumber \\
&&\hspace{-0.5cm}
= \frac{1}{\kappa} \left( \frac{M(\theta_1)}{\sin \theta_2} - \frac{M(\theta_2)}{\sin \theta_1} \right) + \frac{1}{\kappa}  \int_{\theta_1}^{\theta_2}    \frac{M(\theta_1) \, M(\theta_2)}{M(\theta)} \, \frac{\cos \theta} {\sin^2 \theta} \dd  \theta.  
\label{eq:MM}
\end{eqnarray}
The first term is bounded on $\Omega_4$. For the second term we remark that 
$$M(\theta_1) M(\theta_2)/M(\theta) = \exp  \left[ \kappa \left(\frac{1}{\cos \theta} - \left( \frac{1}{\cos \theta_1} + \frac{1}{\cos \theta_2} \right) \right) \right]. $$
Now, $1/\cos \theta$ is a convex function and $\theta \in [\theta_1,\theta_2]$. Assuming, without loss of generality, that $\theta_1 < \theta_2$, we have:
$$ \frac{1}{\cos \theta} \leq \frac{\theta_2 - \theta}{\theta_2 - \theta_1} \, \frac{1}{\cos \theta_1} + \frac{\theta - \theta_1}{\theta_2 - \theta_1} \, \frac{1}{\cos \theta_2}, $$
and we deduce that 
$$\frac{1}{\cos \theta} - \left( \frac{1}{\cos \theta_1} + \frac{1}{\cos \theta_2} \right) \leq \frac{\theta_1 - \theta}{\theta_2 - \theta_1} \, \frac{1}{\cos \theta_1} + \frac{\theta - \theta_2}{\theta_2 - \theta_1} \, \frac{1}{\cos \theta_2} \leq 0.$$
Therefore, $M(\theta_1) M(\theta_2)/M(\theta) \leq 1$ and since $\frac{\cos \theta} {\sin^2 \theta}$ is also bounded for $\theta$ such that $\sin \theta > 1/2$, we deduce that the second term of (\ref{eq:MM}) is also bounded. Consequently, for 
$(\theta_1,\theta_2) \in \Omega_4$, we get 
\begin{equation}
 \Big| \int_{\theta_1}^{\theta_2} \partial_\theta \psi (\theta) \dd  \theta \Big|^2  \, M(\theta_1) \, M(\theta_2) \leq C \int_{\cos \theta >0} | \partial_\theta \psi (\theta) |^2  \, M (\theta) \, \cos^2 \theta \dd  \theta, 
\label{eq:intparthetapsi2}
\end{equation}
which implies, upon integrating with respect to ($\theta_1,\theta_2)$, that 
\begin{equation} 
I_4 \leq C \int_{\cos \theta >0} | \partial_\theta \psi (\theta) |^2  \, M (\theta) \, \cos^2 \theta \dd  \theta . 
\label{eq:intparthetapsi4}
\end{equation}
Now, we take $(\theta_1,\theta_2) \in \Omega_2$ (obviously the case of $(\theta_1,\theta_2) \in \Omega_3$ can be treated similarly). In this case, $I_2$ is decomposed into $I_2=I_2^1 + I_2^2$ where $I_2^1$ and $I_2^2$ correspond to integrals over the domains $ \Omega_2^1 = (- \pi/6, \pi/6) \times (- \pi/2, -\pi/6)$ and $ \Omega_2^2 = (- \pi/6, \pi/6) \times (\pi/6, \pi/2)$ respectively. We consider $I_2^2$ and a similar proof can be made for $I_2^1$. For $(\theta_1,\theta_2) \in \Omega_2^2$, we can decompose:
$$ \int_{\theta_1}^{\theta_2} \partial_\theta \psi (\theta) \dd  \theta = \int_{\theta_1}^{\pi/6} \partial_\theta \psi (\theta) \dd  \theta + \int_{\pi/6}^{\theta_2} \partial_\theta \psi (\theta) \dd  \theta , $$
and consequently we have:
\begin{equation}
 \Big| \int_{\theta_1}^{\theta_2} \partial_\theta \psi (\theta) \dd  \theta \Big|^2 \leq 2 \, \, \Big| \hspace{-1mm} \int_{\theta_1}^{\pi/6} \partial_\theta \psi (\theta) \dd  \theta \Big| ^2+ 2 \, \,  \Big| \hspace{-1mm} \int_{\pi/6}^{\theta_2} \partial_\theta \psi (\theta) \dd  \theta \Big|^2. 
\label{eq:intparthetapsi3}
\end{equation}
Each of the terms will be evaluated separately. For the first term, we proceed as for $I_1$ and notice that $M(\theta)$ and $M(\theta) \, |\cos \theta|$ are bounded from above and from below away from zero, so that we can estimate it similarly to (\ref{eq:inttheta}) and its contribution to $I_2^2$ leads to a similar estimate as (\ref{eq:intparthetapsi}). For the second term, we use (\ref{eq:inttheta2}) with $\theta_1$ replaced by $\pi/6$ and likewise, we find (\ref{eq:inttheta3}) still with $\theta_1$ replaced by $\pi/6$. We can keep on going through all the steps of the proof and we eventually get an estimate similar to (\ref{eq:intparthetapsi2}) with $\theta_1$ replaced by $\pi/6$. This shows that the second term of (\ref{eq:intparthetapsi3}) leads to an estimate similar to (\ref{eq:intparthetapsi4}). Therefore, we get 
\begin{equation} 
I_2 \leq C \int_{\cos \theta >0} | \partial_\theta \psi (\theta) |^2  \, M (\theta) \, \cos^2 \theta \dd  \theta . 
\label{eq:intparthetapsi7}
\end{equation}
Now, inserting (\ref{eq:intparthetapsi}), (\ref{eq:intparthetapsi4}), (\ref{eq:intparthetapsi7}) and its counterpart for $I_3$ into (\ref{eq:poincare3}) leads to  the Poincar\'e estimate (\ref{eq:poincare_2}) and eventually to  (\ref{eq:poincare}), which ends the proof of the Poincar\'e estimate. 

\medskip
\noindent
(iii) We restore the indices $\theta_0$. Then it is obvious that $\psi$ and $\varphi$ satisfy the conditions (\ref{eq:solv}), (\ref{eq:cancel}) respectively and that they belong to ${\mathcal V}_{\theta_0}$. So, they belong to $\tilde {\mathcal V}_{\theta_0}$. Furthermore, by the change of variable $\theta' = \theta + \theta_0$ in the variational formulation satisfied by $\varphi_0$ (i.e.\ with $\theta_0 = 0$), we easily see that $\varphi$ satisfies the variational formulation (\ref{eq:vf1}),  (\ref{eq:vf2}). Since the solution of this variational formulation in the space $\tilde {\mathcal V}_{\theta_0}$ is unique, $\varphi$ is this solution. \endproof

\subsection{Equilibria}
\label{subsec:equi}

We now define the notion of equilibria and characterize them:
\begin{definition}
An equilibrium is a function $f$: $\theta \in {\mathbb R}/(2 \pi {\mathbb Z}) \mapsto f(\theta)$ such that $f \in L^1(-\pi,\pi)$,  $f/M_{\bar \theta_f} \in {\mathcal V}_{\bar \theta_f}$ and $Q^0_\text{al}(f)=0$. We denote by ${\mathcal E}$ the set of equilibria of $Q^0_\text{al}$. 
\label{def:equi}
\end{definition}

\noindent
We also introduce 
\begin{definition}
Let $\rho_+$, $\rho_-$ be real numbers in $[0,\infty)$ and $\bar \theta \in [-\pi/2,\pi/2)$ be an angle of lines. The function $\bar f_{\rho_+,\rho_-, \bar \theta}$ is defined by 
\begin{eqnarray}
\bar f_{\rho_+,\rho_-, \bar \theta}(\theta)&=& \left\{ \begin{array}{lll} 
\rho_+ M_{\bar \theta} (\theta) & \mbox{ for } & \cos (\theta - \bar \theta) >0, \\
\rho_- M_{\bar \theta} (\theta) & \mbox{ for } & \cos (\theta - \bar \theta) <0. \end{array} 
\right.
\label{eq:def_fr+r-}
\end{eqnarray}
It can also be written:
\begin{eqnarray}
\bar f_{\rho_+,\rho_-, \bar \theta}
&=& M_{\bar \theta} (\theta) \, \big[ \rho_+ \, \chi_{\bar \theta}^+ + \rho_- \, \chi_{\bar \theta}^- \big], 
\label{eq:def_fr+r-_2}
\end{eqnarray}
where $\chi_{\bar \theta}^\pm$ are the indicator functions:
\begin{equation}
\chi_{\bar \theta}^\pm (\theta) = \left\{ \begin{array}{ccc} 1 & \mbox{ if } & \pm \cos (\theta - \bar \theta) >0, \\
0 & \mbox{ if } & \pm \cos (\theta - \bar \theta) <0. 
\end{array} \right. 
\label{def:chi}
\end{equation}
\label{def:fr+r-}
\end{definition}

\noindent
Now, with this definition and Prop. \ref{prop:lineareq}, we can state the

\begin{lemma}
\label{lem:vonmises}
The set ${\mathcal E}$ of equilibria of $Q^0_\text{al}$ is given by
$${\mathcal E} =  \big\{ \bar f_{\rho_+,\rho_-, \bar \theta} \, \, | \, \, (\rho_+,\rho_-) \in [0,\infty)^2, \, \, \bar \theta \in [-\pi/2,\pi/2) \big\}. $$ 
\end{lemma}

\medskip 
\noindent 
{\bf Proof.} We use (\ref{eq:defQal}). Thanks to Prop. \ref{prop:lineareq} (i), we first note that, for a given $\theta_0 \in [-\pi/2,\pi/2)$, we have 
\begin{eqnarray*}
 M_{\theta_0} {\mathcal L}_{\theta_0} \Big( \frac{f}{M_{\theta_0}} \Big) = 0 &\Longleftrightarrow&  \frac{f}{M_{\theta_0}} \mbox{ is of the form (\ref{eq:equil})} \\
&\Longleftrightarrow&  \exists  (\rho_+,\rho_-) \in [0,\infty)^2 \mbox{ s.t. } f= \bar f_{\rho_+,\rho_-,\theta_0}.
\end{eqnarray*}
Applying this statement with $\theta_0 = \bar \theta_f$, we deduce that if $f \in {\mathcal E}$, i.e.\ if  $Q^0_\text{al} f =0$, then $f$ is necessarily written as $f= \bar f_{\rho_+,\rho_-,\bar \theta_f}$ and is consequently of the form (\ref{eq:def_fr+r-}). Conversely, for given $(\rho_+,\rho_-) \in [0,\infty)^2$, and $\bar \theta \in [-\pi/2,\pi/2)$, we consider $\bar f_{\rho_+,\rho_-, \bar \theta}$ and need to show that $ Q^0_\text{al} \bar f_{\rho_+,\rho_-, \bar \theta} = 0$. Since $\bar f_{\rho_+,\rho_-, \bar \theta}$ is such that $ M_{\bar \theta} {\mathcal L}_{\bar \theta} ( f/M_{\bar \theta} ) = 0 $ (by the equivalence shown above), it is enough to show that $\bar \theta_{\bar f_{\rho_+,\rho_-, \bar \theta}} = \bar \theta$, or equivalently, that:
$$\int_{-\pi}^{\pi} \sin 2 \big( \theta - \bar \theta \big) \, \bar f_{\rho_+,\rho_-, \bar \theta} \dd  \theta = 0. $$
But we can make the change of variables $\theta' = \theta - \bar \theta$ and use the periodicity to write:
$$\int_{-\pi}^{\pi} \sin 2 \big( \theta - \bar \theta \big) \, \bar f_{\rho_+,\rho_-, \bar \theta} \dd  \theta = \int_{-\pi}^{\pi} \sin 2 \theta \, \bar f_{\rho_+,\rho_-, 0} \dd  \theta,  $$
and this last integral is clearly zero by symmetry. \endproof

\subsection{Generalized Collision Invariants}
\label{subsec_GCI}

We first define the notion of a collision invariant, thanks to the

\begin{definition}
A function $\theta \in {\mathbb R}/(2 \pi {\mathbb Z}) \mapsto \Psi(\theta)$ is a collision invariant (CI) if and only if it satisfies
\begin{align*}
\int Q^0_\text{al}f (\theta) \, \Psi(\theta) \dd  \theta=0 \qquad \text{for all }  f .
\end{align*}
\end{definition}
We deliberately leave the function setting vague as this notion will have to be enlarged later on. Clearly, the set of collision invariants is a vector space. 
At this point, however, we only have a one-dimensional space of collision invariants, spanned by the constant function $\Psi(\theta) = 1$, corresponding to mass conservation. Following the approach in \cite{Degond2008}, we therefore generalize the concept of a collision invariant and define a generalized collision invariant:

\begin{definition}
\label{def:GCI}
Given an angle of lines $\theta_0 \in [-\pi/2,\pi/2)$, a function $\Psi \in {\mathcal V}_{\theta_0}$ is a generalized collision invariant (GCI) associated with $\theta_0$ if and only if
\begin{align}
\langle {\mathcal L}_{\theta_0} \frac{f}{M_{\theta_0}}, \Psi \rangle_{{\mathcal V}_{\theta_0}',{\mathcal V}_{\theta_0}} =0 \quad \text{ for all } \quad  f  \in {\mathcal V}_{\theta_0}\cap L^1(-\pi,\pi) \quad \text{ s.t. } \quad \bar \theta_f= \theta_0. 
\label{eq:defGCI}
\end{align} 
Here the equality $\bar \theta_f= \theta_0$ is an equality of angles of lines i.e. holds modulo $\pi$.  The set of GCI associated with $\theta_0$ is a linear space denoted by ${\mathcal G}_{\theta_0}$. 
\label{def_GCI}
\end{definition}

\noindent
We note that for $f \in {\mathcal V}_{\theta_0}$ such that $ {\mathcal L}_{\theta_0} \frac{f}{M_{\theta_0}} \in  {\mathcal H}_{\theta_0}$, the duality bracket $\langle {\mathcal L}_{\theta_0} \frac{f}{M_{\theta_0}}, \Psi \rangle_{{\mathcal V}_{\theta_0}',{\mathcal V}_{\theta_0}}$ reduces to the classical integral $\int {\mathcal L}_{\theta_0} \frac{f}{M_{\theta_0}}(\theta) \, \Psi(\theta) \, M_{\theta_0}(\theta) \dd  \theta $. Since $\Psi$ being a CI equivalently means that 
\begin{align*}
\int {\mathcal L}_{\bar \theta_f} \frac{f}{M_{\bar \theta_f}}(\theta) \, \Psi(\theta) \, M_{\bar \theta_f}(\theta) \dd  \theta =0 \qquad \text{ for all }  f ,
\end{align*} 
we pass from the definition of a CI to that of a GCI by replacing the orientation $\bar \theta_f$ of $f$ by an arbitrary orientation $\theta_0$ and requesting that the integral is zero only for those functions $f$ whose mean orientation $\bar \theta_f$ coincides with $\theta_0$. Because this means we restrict the set of functions $f$, we expect to find more $\Psi$, fulfilling the definition of a GCI as compared to that of a CI. We can reformulate the definition of a GCI in the following way: 

\begin{lemma}
Given an angle of lines $\theta_0 \in [-\pi/2,\pi/2)$, a function $\Psi \in {\mathcal V}_{\theta_0}$ is a GCI associated with $\theta_0$ if and only if:
\begin{align}
\exists a\in \mathbb{R}  \quad \mbox{ such that } \quad  {\mathcal L}_{\theta_0} \Psi (\theta) = a \sin 2 (\theta - \theta_0) \quad  \mbox{ in } \quad  {\mathcal V}_{\theta_0}'. 
\label{eq:Prob_GCI}
\end{align} 
\label{lem:GCI_equiv}
\end{lemma}

\noindent
{\bf Proof.} We note that the constraint $ \bar \theta_f= \theta_0$ can be equivalently rephrased as the linear constraint: 
$$\int_{-\pi}^{\pi} \sin 2 \big( \theta - \theta_0 \big) \, f(\theta) \dd  \theta = 0. $$
By a classical duality argument already used in \cite{Degond2008}, statement (\ref{eq:defGCI}) is equivalent to the existence of a constant $a \in {\mathbb R}$ such that 
\begin{align*}
\langle {\mathcal L}_{\theta_0} \frac{f}{M_{\theta_0}}, \Psi \rangle_{{\mathcal V}_{\theta_0}',{\mathcal V}_{\theta_0}} = a \int_{-\pi}^{\pi} \sin 2 \big( \theta - \theta_0 \big) \, f(\theta) \dd  \theta \quad \text{ for all }  f  \in {\mathcal V}_{\theta_0}\cap L^1(-\pi,\pi). 
\end{align*} 
Using the formal self-adjointness of $ {\mathcal L}_{\theta_0} $ (see Definition \ref{def:linear_op} (iv)), this is equivalent to
\begin{align*}
\langle {\mathcal L}_{\theta_0} \Psi, \frac{f}{M_{\theta_0}} \rangle_{{\mathcal V}_{\theta_0}',{\mathcal V}_{\theta_0}} = a \Big(
 \sin 2 \big( \theta - \theta_0 \big) ,  \frac{f}{M_{\theta_0}} \Big)_{{\mathcal H}_{\theta_0}} \quad \text{ for all }  f  \in {\mathcal V}_{\theta_0}\cap L^1(-\pi,\pi),
\end{align*} 
which implies that $\Psi$ satisfies (\ref{eq:Prob_GCI}). The converse is obvious. \endproof

We now have the following result:
\begin{proposition}
(i) Given an angle of lines $\theta_0 \in [-\pi/2,\pi/2)$, the space ${\mathcal G}_{\theta_0}$ of GCI associated with $\theta_0$ is of dimension three and is spanned by $\chi_{\theta_0}^+$, $\chi_{\theta_0}^-$ and $g_{\theta_0}$, with 
$\chi_{\theta_0}^\pm$ given by (\ref{def:chi}) and $g_{\theta_0}$, the unique solution in $ {\mathcal V}_{\theta_0}$ of the variational formulation (\ref{eq:vf1}), (\ref{eq:vf2}) with right-hand side $\psi(\theta) = \sin 2 (\theta - \theta_0)$ and satisfying the cancellation condition (\ref{eq:cancel}). \\
(ii) The functions $\chi_{\theta_0}^\pm$ and $g_{\theta_0}$ can be written:
\begin{equation} 
\chi_{\theta_0}^\pm (\theta) = \chi_{0}^\pm (\theta-\theta_0), \qquad g_{\theta_0} (\theta) = g_{0} (\theta-\theta_0).  
\label{eq:transinvar}
\end{equation}
In the remainder, we will write $\chi^\pm$ and $g$ for $\chi^\pm_0$ and $g_0$. In addition, $g$ satisfies the symmetry conditions: 
\begin{equation}
g(-\theta)=-g(\theta), \qquad g(\pi-\theta)=-g(\theta). 
\label{eq:symg}
\end{equation}
\label{prop:GCI}
\end{proposition}

\noindent
{\bf Proof.} (i) First, we note that the solution space of Problem (\ref{eq:Prob_GCI}) is the linear space spanned by the solution of the same problem with $a=1$. Therefore, we restrict ourselves to $a=1$. Now, solving Problem (\ref{eq:Prob_GCI}) with $a=1$ amounts to solving the variational formulation \eqref{eq:vf1}, \eqref{eq:vf2} with right-hand side $\psi(\theta) = \sin 2 (\theta - \theta_0)$. We apply Prop. \ref{prop:lineareq} (ii). A necessary and sufficient condition for the existence of a solution is that the right-hand side satisfies the solvability condition (\ref{eq:solv}), which in the present case, with $\psi = \sin 2 (\theta - \theta_0)$ reads: 
$$\int_{\cos (\theta - \theta_0) >0} \sin 2 (\theta - \theta_0) \, M_{\theta_0} \dd  \theta = 0, \quad  \int_{\cos (\theta - \theta_0) <0} \sin 2 (\theta - \theta_0) \, M_{\theta_0} \dd  \theta = 0.
$$
This condition is obviously satisfied by symmetry. Therefore, there exists a unique solution $g_{\theta_0}$ of the problem satisfying the cancellation condition (\ref{eq:cancel}) and all solutions of the problems differ from this one by the addition of a solution of the homogeneous problem, i.e.\ of the type (\ref{eq:equil}). Such solutions are arbitrary linear combinations of $\chi_{\theta_0}^\pm$. \\
(ii) The statement is obvious for $\chi_{\theta_0}^\pm$. For $g_{\theta_0}$, we simply apply Prop. \ref{prop:lineareq} (iii). Finally, we remark that $\sin 2 \theta$ satisfies the symmetry conditions (\ref{eq:symg}). Therefore, using the uniqueness part of Prop. \ref{prop:lineareq} (ii), we deduce that $g$ satisfies the same symmetry conditions. 
\endproof

\begin{figure}[h!]
\centering
\includegraphics[width=0.5\textwidth]{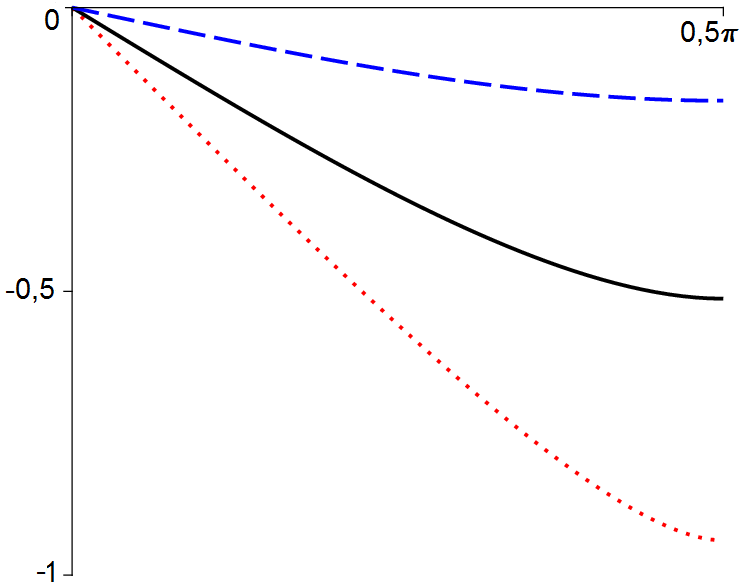}
\caption{$g(\theta)$ for $\kappa=0.5, 2,10$ (red-dotted, black-solid, blue-dashed).}
\label{fig:g}
\end{figure}

\medskip
In the present simple two-dimensional case, we can have an explicit formula for $g$, given by the Lemma below and depicted in Figure \ref{fig:g}. 

\begin{lemma}
For $\theta \in[0,\pi/2]$ we have
\begin{align}
&g(\theta)=-\int_0^\theta \frac{\int_\beta^{\frac{\pi}{2}} \sin 2\alpha \, \exp(-\frac{\kappa}{\cos \alpha }) \dd \alpha}{\cos^2 \beta \, \exp(-\frac{\kappa}{\cos \beta})} \dd  \beta, 
\label{eq:g_explicit}
\end{align} 
and $g$ is extended to ${\mathbb R}/(2 \pi {\mathbb Z})$ by the symmetry conditions (\ref{eq:symg}). 
\label{lem:g_explicit}
\end{lemma}

\noindent
{\bf Proof.} Solving for $g$ consists of solving the following problem:
\begin{align}
\label{eqn:g}
\frac{\dd}{\dd \theta}\Big(\cos^2 \theta \, \,  e^{-\frac{\kappa}{\cos \theta}} \, \frac{\dd g}{\dd\theta} (\theta)\Big)=\sin 2\theta \, \, e^{-\frac{\kappa}{\cos \theta}}. 
\end{align}
Eq. \eqref{eqn:g} can be solved explicitly for $\theta \in[0,\pi/2]$ and then extended using the symmetry conditions (\ref{eq:symg}). Integration constants are determined by the cancellation conditions~(\ref{eq:cancel}). \endproof

\section{Hydrodynamic Model}
\label{sec:hydro}

\subsection{Derivation}
\label{subsec:deriv}

This section is concerned with the limit $\e\rightarrow 0$ in Problem (\ref{eqn:scaling_hydro})-(\ref{eqn:Qal}).  The first step is the observation that the dominant term $Q^0_\text{al} f^\e=0$ requires the limiting distribution to lie in the kernel of $Q^0_\text{al}$, characterized by Lemma \ref{lem:vonmises}. To determine the equations satisfied by the macroscopic quantities $\rho_+(x,t), \rho_-(x,t)$ and $\bar \theta(x,t)$, it is necessary to consider the terms of order $\e$ in (\ref{eqn:scaling_hydro}). This is done by multiplication with the three independent GCI given in Prop. \ref{prop:GCI} and subsequent integration over $\theta$. The computations are similar to those in \cite{Degond2008} and are described in more detail in the Appendix (Section \ref{subsec:app:deriv}). Here, we simply state the result:

\begin{theorem}
\label{thm:limit}
Taking the (formal) limit $\e\rightarrow 0$ in (\ref{eqn:scaling_hydro})-(\ref{eqn:Qal}), we obtain
\begin{align*}
f^\e(x,\theta,t)\, \, \longrightarrow \, \, \bar f_{\rho_+(x,t), \rho_-(x,t), \bar \theta(x,t)}(\theta), 
\end{align*}
where the macroscopic quantities  $\rho_\pm(x,t)$ and $\bar \theta(x,t)$ have values in $[0,\infty)$ and $[-\pi/2,\pi/2)$ respectively and fulfill the following system of equations:
\begin{align}
&\partial_t\rho_+ + d_1 \, \nabla_x\cdot \big(\rho_+ \, v(\bar \theta) \big)=0, \label{eqn:macro_rho+}
\\
&\partial_t\rho_- - d_1 \, \nabla_x\cdot \big(\rho_- v(\bar \theta) \big)=0 , \label{eqn:macro_rho-}\\
&(\rho_+ + \rho_-)\, \partial_t \bar \theta + d_2 \, (\rho_+-\rho_-) \, \big(v(\bar \theta) \cdot \nabla_x\big) \bar \theta + \mu\,  v(\bar \theta)^\perp \cdot \nabla_x(\rho_+-\rho_-) = 0, \label{eqn:macro_theta}
\end{align}
where we recall that $v(\bar \theta) = (\cos \bar \theta, \sin \bar \theta)$ and $v^\bot(\bar \theta) = (-\sin \bar \theta, \cos \bar \theta)$. The coefficients $d_1$, $d_2$ and $\mu$ are given by: 
\begin{align}
&d_1=\langle\cos\rangle_M, \quad d_2=\frac{\big\langle g \frac{\sin}{\cos}\big\rangle_M}{\big\langle g \frac{\sin}{\cos^2}\big\rangle _M}, \quad \mu=\frac{1}{\kappa}\frac{\langle g \sin\rangle _M}{\big\langle g \frac{\sin}{\cos^2}\big\rangle _M}.
\label{eq:def_coef}
\end{align}
where $g$ is the GCI defined in Prop. \ref{prop:GCI} and explicitly given by (\ref{eq:g_explicit}). For a function $\varphi(\theta)$, we denote by $\langle \varphi \rangle_M$ the following average: 
\begin{align}
\label{eqn:average}
\langle \varphi \rangle_M= 2 \int_0^{\pi/2} \varphi(\theta)\, M(\theta) \dd  \theta =\frac{2}{Z_\kappa} \int_0^{\pi/2} \varphi (\theta) \, e^{-\frac{\kappa}{\cos \theta}}\dd  \theta,
\end{align}
where $M$ is $M_{\bar \theta}$ for $\bar \theta = 0$. 
\end{theorem}

\noindent
The proof can be found in the Appendix (Section \ref{subsec:app:deriv}).

\medskip
As shown in the proof of Lemma \ref{lem:vonmises}, the macroscopic quantity $\bar \theta(x,t)$ is the local mean angle of lines as defined in \eqref{eqn:Omega_0} of the equilibrium distribution $\bar f$. The functions $\rho_+(x,t)$ and $\rho_-(x,t)$ describe the local densities of particles that point in direction $[\bar \theta(x,t)-\pi/2, \bar \theta(x,t)+\pi/2)$ and $[\bar \theta(x,t)+\pi/2, \bar \theta(x,t)+3\pi/2)$ respectively. Eqs. (\ref{eqn:macro_rho+}) and (\ref{eqn:macro_rho-}) simply describe number density transport in the direction of $v(\bar \theta)$ (for~$\rho_+$) or~$- v(\bar \theta)$ (for $\rho_-$). These equations express local conservations of the number densities of the particles moving in the direction of $v(\bar \theta)$ or $-v(\bar \theta)$ respectively. Since GCI do not express local conservations (unless they are CI), in general one cannot expect that they will lead to conservation laws. Here, Eqs. (\ref{eqn:macro_rho+}) and (\ref{eqn:macro_rho-}) are conservation laws in spite of the fact that they are derived from using the GCI $\chi_{\bar \theta_f}^\pm$.  This unusual circumstance is due to the particular expression of the collision operator $Q^0_\text{al}$.

Note that whenever $\rho_+=\rho_-$, it follows that $\partial_t \bar \theta =0$ and $\partial_t(\rho_++\rho_-)=0$. An alternative formulation of System (\ref{eqn:macro_rho+})-(\ref{eqn:macro_theta})  can be obtained by defining $\rho=\rho_+ + \rho_-$ and $\delta=\rho_+ - \rho_-$. We get:
\begin{align}
&\partial_t\rho + d_1 \nabla_x\cdot \big(\delta\, v(\bar \theta) \big)=0, \label{eqn:macro_rho_alt}
\\
&\partial_t\delta  -d_1 \nabla_x\cdot \big(\rho\,v(\bar \theta) \big)=0, \label{eqn:macro_delta_alt}\\
&\rho\,\partial_t \bar \theta +d_2\,\delta \, \big(v(\bar \theta) \cdot \nabla_x\big) \bar \theta + \mu \,  \big(v(\bar \theta)^\bot \cdot \nabla_x \big) \delta=0. \label{eqn:macro_theta_alt}
\end{align}

It is instructive to compare with the Self-Organized Hydrodynamic (SOH) system derived in \cite{Degond2008}, which, in dimension two,  takes the form:
 \begin{align}
&\partial_t\rho+\tilde d_1 \nabla_x\cdot \big(\rho\,v(\bar \theta) \big)=0, \label{eq:SOH_rho}\\
&\rho\, \Big( \partial_t \bar \theta +\tilde d_2\, \big(v(\bar \theta) \cdot \nabla_x \big) \bar \theta \Big)+\tilde \mu \big( v(\bar \theta)^\bot \cdot \nabla_x\big) \rho=0, \label{eq:SOH_theta}
\end{align}
where the constants $\tilde d_1, \tilde d_2$ and $\tilde \mu$ depend on $\kappa$ and the collision invariants of the corresponding alignment operator. This model is based on polar alignment. Here by contrast to System (\ref{eqn:macro_rho+})-(\ref{eqn:macro_theta}), $\bar \theta$ is an angle of vectors, i.e. $\bar \theta  \in {\mathbb R}/(2 \pi {\mathbb Z})$ and the particles are distributed according to a classical von Mises-Fisher distribution (i.e.\ their angular distribution is proportional to  $e^{\kappa \cos (\theta-\bar \theta)}$) and they move (with a certain dispersion described by $\kappa$) in the direction of $\bar \theta$.\\

The fact the nematic alignment leads to an equilibrium \emph{angle of lines} as opposed to an \emph{angle of vectors} is reflected in the existence of two opposing local maxima at $\bar \theta$ and $\bar \theta+\pi$ of the GVM $M_{\bar \theta}(\theta)$, whilst the classical von Mises-Fisher distribution has only one. The factor $\cos^2(\theta-\bar \theta_f)$ in front of the diffusion operator in $Q^0_\text{al}$ allows for this separation into two groups and causes $M_{\bar \theta}(\theta)=0$ for $\theta=\bar \theta \pm\pi/2$, i.e. in equilibrium there are no particles moving perpendicular to the equilibrium angle of lines. On the other hand, both the von Mises-Fisher distribution and the GVM share their general "bump-like" shape and the convergence to a Dirac delta for $\kappa \rightarrow \infty$. At the level of the macroscopic equations further similarities and differences become evident. In the nematic case both groups of particles have their own conservation law, consequently we have three instead of two macroscopic equations. The general structure remains the same: Transport equations for the density complemented by an equation describing the angle evolution, consisting of a time derivative, a convection term and a pressure term. However for the nematic SOH model both the convection and the pressure term depend on the local density difference $\delta$. In the case where there is only one group of particles, e.g. $\rho_+=\rho$ and $\rho_-=0$ system \eqref{eqn:macro_rho+}-\eqref{eqn:macro_theta} reduces to \eqref{eq:SOH_rho}-\eqref{eq:SOH_theta} (with different constants), since then the nematic character of the collisions does not play a role.

\subsection{Hyperbolicity of the Macroscopic Model}
\label{subsec:hyperbol}

An important property of the original SOH system (\ref{eq:SOH_rho})-(\ref{eq:SOH_theta}) is its hyperbolicity. In \cite{Frouvelle2012} a variant of the SOH model was examined. It was found that in some regions hyperbolicity is lost, which affects the behavior of the solutions. Therefore, we examine the hyperbolicity of the macroscopic model (\ref{eqn:macro_rho+})-(\ref{eqn:macro_theta}). To simplify the analysis we rescale the system with $t \rightarrow d_1 t$. In the following we therefore abbreviate $\hat d_2=d_2 /d_1$ and $\hat \mu=\mu /d_1$.

\begin{lemma}
The macroscopic system  (\ref{eqn:macro_rho+})-(\ref{eqn:macro_theta}) is hyperbolic.
\end{lemma}
\medskip \noindent {\bf Proof.}
We use the alternative formulation (\ref{eqn:macro_rho_alt})-(\ref{eqn:macro_theta_alt}), which leads to simpler formulas. 
we define $u(x,t) = (\rho(x,t), \delta(x,t), \bar \theta(x,t))$ with $x=(x_1,x_2)\in \mathbb{R}^2$, yielding
\begin{align*}
&\partial_t u+A \partial_{x_1} u+B \partial_{x_2} u=0, 
\end{align*}
where
\begin{align*}
&A(\rho, \delta,\bar \theta)=\begin{pmatrix}
0 && \cos \bar \theta &&  -\delta \, \sin \bar \theta\\
\cos \bar \theta && 0 && -\rho \, \sin \bar \theta \\
0 && -\frac{\hat \mu}{\rho}\,\sin\bar \theta && \frac{\hat d_2 \delta}{\rho} \, \cos \bar \theta 
\end{pmatrix}, \quad
B(\rho, \delta,\bar \theta)=A(\rho, \delta,\bar \theta-\pi/2).
\end{align*}
To determine the hyperbolicity of the system, it suffices to show that the matrix $A$ has real and distinct eigenvalues for all $\bar \theta \in[-\pi/2, \pi/2)$. The characteristic polynomial of $A$, $\text{charPoly}_A(z)$ is of third degree and has the shape
\begin{align*}
\text{charPoly}_A(z)=z^3-\frac{\delta}{\rho}\cos{\theta}\,\hat d_2\, z^2-\left(\hat \mu\sin^2 \theta+\cos^2 \theta\right)z-\frac{\delta}{\rho}\cos \theta\left(\mu \sin^2\theta -\hat d_2 \cos^2 \theta\right).
\end{align*}
To understand its roots, we use the following definition of a discriminant of a third order polynomial: For $p(z)=a_3 z^3+a_2 z^2+a_1 z+a_0$ we define
\begin{align*}
\text{discr}_p=-27 a_3^2 a_0^2+18 a_3 a_0 a_1 a_2+a_2^2 a_1^2-4 a_2^3 a_0-4 a_1^3 a_3.
\end{align*}
Next we observe that the discriminant of the characteristic polynomial of $A$, which we define as $\Delta=\text{discr}_{\text{charPoly}_A}$, can be interpreted as a quadratic polynomial in $X=(\delta/\rho)^2$ with coefficients dependent on $c=\cos^2 \theta$, namely:
\begin{align*}
\Delta(X)&=\alpha(c) X^2+\beta(c) X+\gamma(c), 
\end{align*}
where
\begin{align*}
\alpha(c)&=4\,\hat d_2^3\, c^2 \left[c\,(\hat d_2+\hat \mu)-\hat \mu\right], \\
\beta(c)&= c\left[\left(\hat d_2^2\hat \mu^2-18\,\hat d_2\hat \mu^2-27\hat \mu^2-36\,\hat d_2\hat \mu-8\,\hat d_2^2-20\,\hat d_2^2\hat \mu\right)c^2\right.\\
&+\left.2\hat \mu\left(-\hat d_2^2\hat \mu+18\,d_2+10\,\hat d_2^2+27\hat \mu+18\hat \mu\,\hat d_2\right)c-\hat \mu^2(18\,\hat d_2-\hat d_2^2+27)\right], \\
\gamma(c)&=4\left[c\, (1-\hat \mu)+\hat \mu\right]^3.
\end{align*}
If $\Delta(X)>0$ for all $c\in[0,1]$ and all $X \in[0,1]$ we have three distinct, real eigenvalues and the system is (strictly) hyperbolic. The roots of $\Delta(X)$ are given by 
$$X_{1,2}(c)=\frac{-\beta(c)\pm \sqrt{\beta(c)^2-4 \alpha(c)\gamma(c)}}{2\alpha(c)}.$$ 
Simple calculations show that for $c>\frac{9(\hat \mu+\hat d_2)}{9 \hat \mu+\hat \mu \hat d_2+8 \hat d_2}$ the discriminant $\Delta(X)$ has no real roots and is positive for all $X\in[0,1]$. Otherwise we can distinguish between $c>s:=\frac{\hat \mu}{\hat \mu+\hat d_2}$ and $0\leq c\leq s$. In the first case both roots $X_{1,2}(c)$ are negative are therefore not relevant. In the second case we lose hyperbolicity for $X>X_2(c)$. However a simple analysis shows that $X_2(c)\geq 1$ for all $0\leq c\leq s$. In fact $X_2(c)\rightarrow \infty$ for $c\rightarrow 0^+$ and $c\rightarrow s^-$ and $\min_{c\in[0,s]}X_2(c)=1$. Since $X\in[0,1]$ this shows the hyperbolicity of the system. 
\endproof

\section{Weakly non-local scaling}
\label{sec:nonlocal}

In this section we want to investigate how the model changes, if we allow for a relatively larger interaction radius. With the scaling choice below, this leads to interactions that are no longer local in the limit, but rather weakly non-local. This and other scaling choices for different models of self-propelled particles have been examined in \cite{Degond2013}.
In the scaling performed in Section \ref{ssec:meanfield}, we used $\hat R=\mathcal{O}(\e)$. If we use $\hat R=\mathcal{O}(\sqrt{\e})$ instead, namely $\hat R = \sqrt{\varepsilon} \, r$,  we are led to the following expression for $J^\varepsilon_f$: 
\begin{eqnarray}
&&\hspace{-1cm}
J^\varepsilon_f(x,t)= \int_{-\pi}^{\pi} \int_{|y-x|\leq \sqrt{\varepsilon} \, r} v( 2 \theta ) \, f(y,\theta,t) \dd y \dd  \theta, 
\label{eqn:scaling_jf_3}
\end{eqnarray}
and $\bar \Theta^\varepsilon_f$ is still given by (\ref{eqn:Thetaf_scaled}). 
Then, we need to consider higher order terms in the Taylor expansion of $J^\varepsilon_f$ and $\bar \Theta^\varepsilon_f$. These are given by the following:

\begin{lemma}
We have 
\begin{eqnarray}
&&\hspace{-1cm}
J^\varepsilon_f(x,t)= \varepsilon \pi r^2 \, \big(  j_f(x,t) + \varepsilon \,  j_f^1(x,t) + {\mathcal O}(\varepsilon^2) \big) , \label{eq:expan_jeps_1storder}\\
&&\hspace{-1cm}
 j_f^1(x,t) = k \, \Delta_x j_f(x,t), \qquad  k = \frac{1}{4 \pi r^2} \int_{|\xi|\leq r} |\xi|^2 \dd  \xi = \frac{r^2}{8}, \label{eq:j1} \\
&&\hspace{-1cm}
\bar \Theta_f^\varepsilon = \bar \theta_f + \varepsilon \bar \theta_f^1 + {\mathcal O}(\varepsilon^2) , \quad  \bar \theta_f^1= \frac{1}{2|j_f|} \big(j_f^1 \cdot v^\bot(2 \bar \theta_f) \big), \label{eq:theta1} 
\end{eqnarray}
and where $j_f$, $\bar \theta_f$ are  given by (\ref{eqn:Omega}) and (\ref{eqn:Omega_0}) respectively.
\label{lem:expansion_bartheta}
\end{lemma}

\noindent
{\bf Proof.}
We change variables to $y = x + \sqrt{\varepsilon} \xi$ and get 
\begin{eqnarray*}
&&\hspace{-1cm}
J^\varepsilon_f(x,t)= \varepsilon \int_{|\xi|\leq r} \int_{-\pi}^{\pi} v( 2 \theta )  \, f(x + \sqrt{\varepsilon} \xi ,\theta,t) \dd \xi \dd  \theta. 
\end{eqnarray*}
By Taylor expansion of $f$ with respect to $\sqrt{\varepsilon}$ and using cancellations of the order $\sqrt{\varepsilon}$ terms due to symmetry, we 
find (\ref{eq:expan_jeps_1storder}), (\ref{eq:j1}). Therefore, we get
\begin{eqnarray}
&&\hspace{-1cm}
\frac{J_f}{|J_f|}=\frac{j_f}{|j_f|} + \varepsilon \frac{1}{|j_f|} \Big(j_f^1 - \big( j_f^1 \cdot \frac{j_f}{|j_f|} \big) \frac{j_f}{|j_f|} \Big)+ {\mathcal O}(\varepsilon^2) \nonumber \\
&&\hspace{-0.5cm}
= v(2 \bar \theta_f) + \varepsilon \frac{1}{|j_f|} \big(j_f^1 \cdot v^\bot(2 \bar \theta_f) \big) \, v^\bot(2 \bar \theta_f) + {\mathcal O}(\varepsilon^2) . \label{eq:auxi_expan_1}
\end{eqnarray}
Expanding $\bar \Theta_f$ according to the first equation of (\ref{eq:theta1}), we get
\begin{eqnarray}
&&\hspace{-1cm}
\frac{J_f}{|J_f|}= v \big( 2 \bar \Theta_f \big) = v \big( 2 \bar \theta_f \big)  + 
2\, \varepsilon\, \bar \theta_f^1 \, v^\bot\big( 2 \bar \theta_f \big) + {\mathcal O}(\varepsilon^2) . 
\label{eq:auxi_expan_2}
\end{eqnarray}
Identifying the coefficients of $\varepsilon$ in (\ref{eq:auxi_expan_1}) and in (\ref{eq:auxi_expan_2}), we get the second equation of (\ref{eq:theta1}).
\endproof

\medskip
From the expansion given in Lemma \ref{lem:expansion_bartheta}, we deduce an expansion of the operators $\tilde Q^\varepsilon_\text{al}f$. We recall that their expressions are given by (\ref{eqn:Qal_scaled}). This expansion is given in the 

\begin{lemma}
We have 
\begin{eqnarray}
&&\hspace{-1cm}
\tilde Q^\varepsilon_\text{al}f = Q^0_\text{al}f + \varepsilon Q^1_\text{al}f + {\mathcal O}(\varepsilon^2), 
\label{eq:expan_Qal}
\label{eq:expan_Qrev}
\end{eqnarray}
where $Q^0_\text{al}f$ is given by (\ref{eqn:Qal}) and $Q^1_\text{al}f$ is given by 
\begin{eqnarray}
&&\hspace{-2cm}
Q^1_\text{al}f=  \bar \theta_f^1 \, \partial_\theta  \bigg[ \Big( 2 \, \delta \big(\cos(\theta- \bar \theta_f) \big) \, \sin^2 (\theta- \bar \theta_f) 
- \big|\cos(\theta- \bar \theta_f) \big| \, \Big)  f \nonumber \\
&&\hspace{6.cm}
+ 2D \,  \cos(\theta- \bar \theta_f) \, \sin(\theta- \bar \theta_f) \, \partial_\theta f \bigg],
\label{eq:Qal1}
\end{eqnarray}
where  $ \delta \big(\cos(\theta- \bar \theta_f) \big) $ is the Dirac delta at $\cos(\theta- \bar \theta_f) = 0$. 
\label{lem:expan_Qal_Qrev}
\end{lemma}

\noindent
{\bf Proof.} Formula (\ref{eq:expan_Qal}) follows from the first equation of (\ref{eq:theta1}) as well as 
$$ \theta - \bar \Theta^\varepsilon_f = \theta - \bar \theta_f - \varepsilon \bar \theta_f^1+ {\mathcal O}(\varepsilon^2).$$
Then, inserting this expansion in (\ref{eqn:Qal_scaled}) and using elementary calculus, we are led to (\ref{eq:Qal1}). Formula (\ref{eq:expan_Qrev}) is obvious. \endproof

\noindent
From the previous two lemmas, in the new scaling, the mean-field kinetic equation (\ref{eqn_meanfield_scaled}) is written: 
\begin{eqnarray}
&&\hspace{-1cm}
\varepsilon \Big( \partial_t f^\varepsilon+ \nabla_x\cdot \big(v(\theta) \, f^\varepsilon \big) - Q^1_\text{al}f^\varepsilon \Big)= Q^0_\text{al}f^\varepsilon + {\mathcal O}(\varepsilon^2), \label{eqn_meanfield_scaled_2} 
\end{eqnarray}
and we will now drop the ${\mathcal O}(\varepsilon^2)$ term. The main theorem of this section is the following.

\begin{theorem}
\label{thm:nonlocal_limit}
Taking the (formal) limit $\e\rightarrow 0$ in \eqref{eqn_meanfield_scaled_2}, we obtain
\begin{align*}
f^\e(x,\theta,t)\rightarrow \bar f_{\rho_+(x,t), \rho_-(x,t), \bar \theta(x,t)}(\theta), 
\end{align*}
with $\bar f_{\rho_+, \rho_-, \bar \theta}$ given by (\ref{eq:def_fr+r-}). The macroscopic quantities $\rho_\pm$ and $\bar \theta$ fulfill
\begin{eqnarray}
&&\hspace{-1.5cm}
\partial_t\rho_+ + d_1 \, \nabla_x\cdot \big(\rho_+ \, v(\bar \theta) \big)=0 , \label{eqn:macro_diff_rho+} \\
&&\hspace{-1.5cm}
\partial_t\rho_- - d_1 \, \nabla_x\cdot \big(\rho_- v(\bar \theta) \big)=0 , \label{eqn:macro_diff_rho-}\\
&&\hspace{-1.5cm}
(\rho_+ + \rho_-)\, \partial_t \bar \theta + d_2 \, (\rho_+-\rho_-) \, \big(v(\bar \theta) \cdot \nabla_x\big) \bar \theta + \mu\,  v(\bar \theta)^\perp \cdot \nabla_x(\rho_+-\rho_-) \nonumber \\
&&\hspace{3cm} 
= 2k \, \Big( d_2 + d_3 -  \mu - \frac{2}{\kappa}   \Big) \, 
\nabla_x \cdot \big( (\rho_+ + \rho_-) \nabla_x \bar \theta \big), 
\label{eqn:macro_nonlocal}
\end{eqnarray}
where $d_1$, $d_2$ and $\mu$ are given by (\ref{eq:def_coef}) and 
\begin{align}
&d_3=2\frac{\big\langle g\frac{\sin^3}{\cos^3}\big\rangle_M}{\big\langle g\frac{\sin}{\cos^2}\big\rangle_M}. \label{eqn:d3}
\end{align}
The definition of $\langle .\rangle_M$ can be found in \eqref{eqn:average}.
\end{theorem}

\noindent
The proof of this theorem is given in the Appendix (see Section \ref{app:nonlocal}). 
The weak non-locality of the interactions leads to the addition on a diffusion term on the right-hand side compared to purely local interactions. The sign of the diffusion constant depends on the sign of $\mathcal{D}=d_2 + d_3 -  \mu - \frac{2}{\kappa}$. Lemma \ref{lem:diffusion_constant_positive} shows that $\mathcal{D}$ is positive.

\begin{lemma}
\label{lem:diffusion_constant_positive}
Let the diffusion constant $\mathcal{D}$ be defined by
\begin{align*}
\mathcal{D}:=d_2+d_3-\mu-\frac{2}{\kappa},
\end{align*}
where $d_2$ and $\mu$ are given by \eqref{eq:def_coef} and $d_3$ given  by \eqref{eqn:d3}. Then it holds that
\begin{align*}
\mathcal{D}\geq 0\quad \text{for}\quad  \kappa\geq 0.
\end{align*}
\end{lemma}
{\bf Proof.}
We start by rewriting
\begin{align}
\mathcal{D}&=d_2+d_3-\mu-\frac{2}{\kappa}\nonumber\\
&=\frac{\big\langle g \frac{\sin}{\cos}\big\rangle_M}{\big\langle g \frac{\sin}{\cos^2}\big\rangle_M}+2\,\frac{\big\langle g \frac{\sin^3}{\cos^3}\big\rangle_M}{\big\langle g\frac{\sin}{\cos^2}\big\rangle_M}-\frac{1}{\kappa}\frac{\big\langle g\sin\big\rangle_M}{\big\langle g\frac{\sin}{\cos^2}\big\rangle_M}-\frac{2}{\kappa}\nonumber\\
&=\frac{1}{\kappa \big\langle g \frac{\sin}{\cos^2}\big\rangle_M } \Bigg\langle g \left[ \kappa \frac{\sin}{\cos}+2\kappa \frac{\sin^3}{\cos^3} -\sin-2\frac{\sin}{\cos^2}\right]\Bigg\rangle_M=\frac{A}{B}.\label{eqn:proof_pos_AB}
\end{align}
Next, using the definition of $\langle . \rangle_M$ given in \eqref{eqn:average}, we examine the numerator $A$. Note that since $Z_\kappa$ would appear in both $A$ and $B$, we can omit it.
\begin{align}
A&=\bigg\langle g \frac{\sin}{\cos^3} \left[ \kappa \left(\cos^2+2\sin^2\right) -\cos \left(\cos^2+2\right)\right]\bigg\rangle_M\nonumber\\
&=\int_0^{\pi/2}\!\!\!  g(\theta) \frac{\sin\theta}{\cos^3\theta}\, e^{-\frac{\kappa}{\cos\theta}} \left[ \kappa \left(\cos^2\theta+2\sin^2\theta\right) -\cos\theta \left(\cos^2\theta+2\right)\right]\dd \theta\nonumber\\
&=\int_0^1 \frac{G(y)}{y^3} e^{-\frac{\kappa}{y}} \left[ \kappa \left(2-y^2\right) -y \left(y^2+2\right)\right] \dd y.\label{eqn:proof_pos_1}
\end{align}
In the last line we used the change of variable $y=\cos\theta$ and the definition
\begin{align*}
G(y)=g(\arccos(y)).
\end{align*}
The last line \eqref{eqn:proof_pos_1} can be split into
\begin{align}
A&=\int_0^1 \frac{G(y)}{y} \frac{\dd}{\dd y} \left(e^{-\frac{\kappa}{y}}\right) \left(2-y^2\right) \dd y-\int_0^1 \frac{G(y)}{y^2} e^{-\frac{\kappa}{y}} \left(y^2+2\right) \dd y\nonumber\\
&=A_1+A_2.\label{eqn:proof_pos_A}
\end{align}
We start by rewriting $A_1$ using integration by parts. This yields
\begin{align*}
A_1&=\frac{G(y)}{y}e^{-\frac{\kappa}{y}} \left(2-y^2\right)\Big|_0^1-\int_0^1 e^{-\frac{\kappa}{y}} \frac{\dd}{\dd y} \left(G(y)\frac{2-y^2}{y}\right)\dd y\\
&=-\int_0^1 e^{-\frac{\kappa}{y}} \frac{2-y^2}{y}\frac{\dd}{\dd y} G(y)\dd y-\int_0^1 e^{-\frac{\kappa}{y}} G(y)\frac{\dd}{\dd y} \left(\frac{2-y^2}{y}\right)\dd y\\
&=A_1^1+A_1^2,
\end{align*}
where the second equality is obtained by noting that the boundary term vanishes. Note that since $ \frac{\dd g}{\dd \theta}\leq 0$ on $[0, \pi/2]$, which can be seen from its definition in \eqref{eq:g_explicit}, it also holds that
\begin{align}
\frac{\dd}{\dd y} G(y)=-\frac{1}{\sqrt{1-y^2}} \frac{\dd g}{\dd \theta}(\arccos(y))\geq 0 \quad \text{for}\quad y\in [0,1].\label{eqn:proof_pos_dG}
\end{align}
We continue by expanding $A_1^2$:
\begin{align*}
A_1^2=\int_0^1 e^{-\frac{\kappa}{y}} G(y)\left(\frac{2+y^2}{y^2}\right)\dd y,
\end{align*}
which shows that in fact $A_1^2=-A_2$ as defined in \eqref{eqn:proof_pos_A}. Using this and \eqref{eqn:proof_pos_dG} we can conclude
\begin{align*}
A=A_1^1&=-\int_0^1 e^{-\frac{\kappa}{y}} \frac{2-y^2}{y}\frac{\dd}{\dd y} G(y)\dd y\leq 0.
\end{align*}
Finally, since $g\leq 0$ on $[0, \pi/2]$, we observe that the denominator $B$ as defined in \eqref{eqn:proof_pos_AB} (again omitting $Z_\kappa$) fulfills
\begin{align*}
B=\kappa \Big\langle g \frac{\sin}{\cos^2} \Big\rangle_M=\kappa \int_0^{\pi/2} g(\theta) e^{-\frac{\kappa}{\cos\theta}}\frac{\sin\theta}{\cos^2 \theta}\dd \theta \leq 0, 
\end{align*}
which shows the claim and finishes the proof.
\endproof


\section{Example: Myxobacteria}
\label{sec:myxo}

In this section we want to give an example of how to use and adapt the derived nematic SOH model to a particular biological problem, the movement of myxobacteria. This type of rod-shaped bacteria can glide without the aid of flagella on surfaces. Under starvation conditions, prior to the formation of fruiting bodies, they go through a so-called \emph{rippling phase}. During this stage traveling waves of high density bands are observed that seemingly pass through each other unaffected. However upon close examination it has been found that a large number of bacteria at the interface of two colliding waves in fact reverse their direction. These reversals can happen spontaneously, but a higher number of head-to-head contacts with other, oppositely moving myxobacteria increases their likelihood \cite{Igoshin2004, Welch2001}. Biologically these reversals are quite fascinating, because they happen by the reorganization of the internal movement machinery as opposed to an actual turning process. The interaction of two bacteria that don't reverse can be described by nematic alignment. In terms of other communication between bacteria no diffusible chemotactic signal is known to coordinate their movement. At this point the nematic SOH model can describe the (nematic) alignment of two colliding bacteria, however the reversals still need to be included.

\subsection{Particle and Mean Field Model}

\noindent
\paragraph{Particle Model}\\
 We complete the description of the particle model from Subsection \ref{subsec:particle} by assuming that particles reverse their direction of motion with a frequency, dependent on the density of oppositely moving bacteria $\rho_\pm^i$. So we let $\Theta_i$ be changed into $\Theta_i+\pi$ at Poisson random times with frequency  $\lambda(\rho_{-\sgn  \cos (\Theta_i - \bar \Theta_i)}^i)$, where 
$$ 
\rho_\pm^i:=\frac{1}{R^2 \pi} \, \, \mbox{Card} \{ k\in\mathcal{N}\, \, \big| \,  \, |X_k-X_i|\leq R \mbox{ and } \pm \cos (\Theta_k - \bar \Theta_i) \geq 0 \}, 
$$
where Card stands for the cardinal of a set. In other words, we have 
\begin{eqnarray}
&&\hspace{-1cm}
\mbox{Prob} \big\{ \Theta_i(t+ T) = \Theta_i(t) + \pi \, \, | \, \, \Theta_i(t+ \tau) = \Theta_i(t) \,\, \forall \tau \mbox{ with } 0 \leq \tau < T \big\}  \nonumber \\
&&\hspace{6cm}
= 1 - \exp \big( - \lambda(\rho_{-\sgn  \cos (\Theta_i - \bar \Theta_i)}^i) \, T \big).
\label{eqn:particle_3}
\end{eqnarray}
(assuming that $\rho_{-\sgn  \cos (\Theta_i - \bar \Theta_i)}^i= \rho_{-\sgn  \cos (\Theta_i - \bar \Theta_i)}^i(t)$ stays constant in the time interval $[t,t+T]$). In these formula, the function $\lambda(\rho_\pm)$ is a reversal frequency. The way it depends on $\rho_\pm$ will be discussed later. \\

\noindent
\paragraph{Mean Field Model}\\
The above considerations change equation \eqref{eqn_meanfield} to 
\begin{equation}
\partial_t f+ v_0\nabla_x\cdot (v(\theta) \, f)=\tilde Q_\text{al}f+\tilde Q_\text{rev}f, \label{eqn_meanfield_myxo}
\end{equation}
where the collision operators $\tilde Q_\text{al}$ is given by \eqref{eqn:scaling_Qal} and 
\begin{eqnarray}
\tilde Q_\text{rev}f \, (x,\theta,t) =\begin{cases}
\lambda(\sigma_{f,+}) \, f(x,\theta+ \pi, t) - \lambda(\sigma_{f,-}) \, f(x,\theta,t),  \, \, \, \text{for} \, \cos(\theta-\bar \Theta_f) > 0, 
\\
\lambda(\sigma_{f,-}) \, f(x,\theta+ \pi, t) -\lambda(\sigma_{f,+}) \,  f(x,\theta, t),  \, \, \, \text{for} \, \cos(\theta-\bar \Theta_f) < 0,  \\
\end{cases} \label{eqn:scaling_Qrev}
\end{eqnarray}
and where $ \bar \Theta_f$ is given by \eqref{eqn:scaling_jf_0}-\eqref{eqn:scaling_jf} and $\sigma_{f,\pm}$ given by
\begin{eqnarray}
\sigma_{f,\pm}(x,t)=\frac{1}{R^2\pi}\int_{\pm \cos(\theta - \bar \Theta_f) \geq 0} \int_{|x-y|\leq R} f(y,\theta,t) \dd y \dd \theta. \label{eqn:scaling_sigmaf}
\end{eqnarray}
To obtain local interactions, we assume as in Section \ref{ssec:meanfield} an interaction radius fulfilling $\hat R=\mathcal{O}(\e)$ and choose for the reversal frequency the scaling $\hat \lambda(\sigma_{\hat f, \pm})=t_0' \lambda(\sigma_{f,\pm})=\frac{t_0}{\e}\lambda(\sigma_{f,\pm})$. The local density is then given by 
\begin{align*}
\sigma^\varepsilon_{f,\pm}(x,t)=\frac{1}{\varepsilon^2 r^2\pi}\int_{\pm \cos(\theta - \bar \Theta^\varepsilon_f) > 0} \int_{|x-y|\leq \varepsilon \, r} f(y,\theta,t) \dd y \dd \theta, 
\end{align*}
and it can be shown that
\begin{align*}
\sigma^\varepsilon_{f,\pm}(x,t)=\rho_{f,\pm}(x,t)+\mathcal{O}(\e^2),
\end{align*}
where
\begin{align}
\rho_{f,\pm}(x,t)=\int_{\pm \cos(\theta - \bar \theta_f) \geq 0}  f(x,\theta,t)  \dd \theta. \label{eqn:rho}
\end{align}
This changes equation \eqref{eqn:scaling_hydro} into 
\begin{eqnarray}
&&\hspace{-1cm}
\e\big( \partial_t f^\varepsilon+\nabla_x \cdot (v(\theta) \,f^\varepsilon) \big)=Q^0_\text{al} f^\varepsilon+\e\, Q^0_\text{rev} f^\varepsilon, \label{eqn:scaling_hydro_myxo}
\end{eqnarray}
where $Q^0_\text{al} f$ is given by \eqref{eqn:Qal} and 

\begin{align}
Q^0_\text{rev}f(x,\theta,t) =\begin{cases}
\lambda(\rho_{f,+}) \, f(x,\theta+ \pi, t) - \lambda(\rho_{f,-}) \, f(x,\theta,t) \, \,  \, \text{for} \,  \cos(\theta-\bar \theta_f) > 0, 
\\
\lambda(\rho_{f,-}) \, f(x,\theta+ \pi, t) -\lambda(\rho_{f,+}) \, f(x,\theta, t) \, \,  \, \text{for} \,  \cos(\theta-\bar \theta_f) < 0, \\
\end{cases} 
\label{eqn:Qrev}
\end{align}
where $\bar \theta_f$ and $\rho_{f,\pm}$ are given by \eqref{eqn:Omega_0} and \eqref{eqn:rho} respectively.

\begin{remark} In \eqref{eqn:scaling_hydro_myxo} the contributions from $Q^0_\text{rev}$ appear at the same scale as the transport operator, which is a consequence of the scaling choice for the reversal rate $\lambda$. Using $\hat \lambda(\sigma_{f,\pm})=t_0 \lambda(\sigma_{f,\pm})$ instead, would lead to a different model, in which $Q^0_\text{al}$ and $Q^0_\text{rev}$ would be of the same order. Consequently also the shape of the equilibria and the macroscopic model would look very different. We leave this path for future work.
\end{remark}

\subsection{Hydrodynamic Limits}
\noindent
\paragraph{Macroscopic Model - Local Scaling}

\begin{proposition}
\label{thm:limit_myxo}
Taking the (formal) limit $\e\rightarrow 0$ in (\ref{eqn:scaling_hydro_myxo})-(\ref{eqn:Qrev}), we obtain
\begin{align*}
f^\e(x,\theta,t)\, \, \longrightarrow \, \, \bar f_{\rho_+(x,t), \rho_-(x,t), \bar \theta(x,t)}(\theta), 
\end{align*}
where the macroscopic quantities  $\rho_\pm(x,t)$ and $\bar \theta(x,t)$ have values in $[0,\infty)$ and $[-\pi/2,\pi/2)$ respectively and fulfill the following system of equations:
\begin{align}
&\partial_t\rho_+ + d_1 \, \nabla_x\cdot \big(\rho_+ \, v(\bar \theta) \big)=\lambda(\rho_+) \, \rho_- - \lambda(\rho_-)\, \rho_+ , \label{eqn:macro_rho+_myxo}
\\
&\partial_t\rho_- - d_1 \, \nabla_x\cdot \big(\rho_- v(\bar \theta) \big)=\lambda(\rho_-)\, \rho_+ - \lambda(\rho_+)\, \rho_- , \label{eqn:macro_rho-_myxo}\\
&(\rho_+ + \rho_-)\, \partial_t \bar \theta + d_2 \, (\rho_+-\rho_-) \, \big(v(\bar \theta) \cdot \nabla_x\big) \bar \theta + \mu\,  v(\bar \theta)^\perp \cdot \nabla_x(\rho_+-\rho_-) = 0, \label{eqn:macro_theta_myxo}
\end{align}
where $v(\bar \theta) = (\cos \bar \theta, \sin \bar \theta)$ and $v^\bot(\bar \theta) = (-\sin \bar \theta, \cos \bar \theta)$. The coefficients $d_1$, $d_2$ and $\mu$ are defined in \eqref{eq:def_coef}.
\end{proposition}

The proof consists of some simple additions to the proof of Theorem \ref{thm:limit} and is given in the Appendix (see Section \ref{ssec:proof_thm:limit_myxo}).\\

\noindent
\paragraph{Macroscopic Model - Weakly non-local Scaling}\\
To obtain the hydrodynamic limit of the myxobacteria model for the weakly non-local scaling described in Section \ref{sec:nonlocal}, we note that in this scaling
\begin{align*}
\sigma^\varepsilon_{f,\pm}(x,t)=\frac{1}{\varepsilon r^2\pi}\int_{\pm \cos(\theta - \bar \Theta^\varepsilon_f) > 0} \int_{|x-y|\leq \sqrt{\varepsilon} \, r} f(y,\theta,t) \dd y \dd \theta. 
\end{align*}
It holds that
\begin{align*}
\sigma^\varepsilon_{f,\pm}(x,t)=\rho_{f,\pm}(x,t)+\mathcal{O}(\e).
\end{align*}
We observe that \eqref{eqn_meanfield_scaled_2} is changed into 
\begin{align}
\varepsilon \Big( \partial_t f^\varepsilon+ \nabla_x\cdot \big(v(\theta) \, f^\varepsilon \big) - Q^1_\text{al}f^\varepsilon-Q^0_\text{rev} \Big)= Q^0_\text{al}f^\varepsilon + {\mathcal O}(\varepsilon^2), \label{eqn_meanfield_scaled_2_myxo} 
\end{align}
with $Q^1_\text{al}$ and $Q^0_\text{rev}$ given by \eqref{eq:Qal1} and \eqref{eqn:Qrev} respectively. The resulting limit is summarized in 

\begin{proposition}
\label{thm:nonlocal_limit_myxo}
Taking the (formal) limit $\e\rightarrow 0$ in (\ref{eqn:scaling_hydro_myxo})-(\ref{eqn:Qrev}), we obtain
\begin{align*}
f^\e(x,\theta,t)\rightarrow \bar f_{\rho_+(x,t), \rho_-(x,t), \bar \theta(x,t)}(\theta), 
\end{align*}
with $\bar f_{\rho_+, \rho_-, \bar \theta}$ given by (\ref{eq:def_fr+r-}). The macroscopic quantities $\rho_\pm$ and $\bar \theta$ fulfill
\begin{eqnarray*}
&&\hspace{-1.5cm}
\partial_t\rho_+ + d_1 \, \nabla_x\cdot \big(\rho_+ \, v(\bar \theta) \big)=\lambda(\rho_+) \, \rho_- - \lambda(\rho_-)\, \rho_+ , \label{eqn:macro_diff_rho+_myxo} \\
&&\hspace{-1.5cm}
\partial_t\rho_- - d_1 \, \nabla_x\cdot \big(\rho_- v(\bar \theta) \big)=\lambda(\rho_-)\, \rho_+ - \lambda(\rho_+)\, \rho_- , \label{eqn:macro_diff_rho-}\\
&&\hspace{-1.5cm}
(\rho_+ + \rho_-)\, \partial_t \bar \theta + d_2 \, (\rho_+-\rho_-) \, \big(v(\bar \theta) \cdot \nabla_x\big) \bar \theta + \mu\,  v(\bar \theta)^\perp \cdot \nabla_x(\rho_+-\rho_-) \nonumber \\
&&\hspace{3cm} 
= 2k \, \Big( d_2 + d_3 -  \mu - \frac{2}{\kappa}   \Big) \, 
\nabla_x \cdot \big( (\rho_+ + \rho_-) \nabla_x \bar \theta \big), 
\label{eqn:macro_nonlocal_myxo}
\end{eqnarray*}
where $d_1$, $d_2$ and $\mu$ are given by (\ref{eq:def_coef}) and $d_3$ is given by \eqref{eqn:d3}. The definition of $\langle .\rangle_M$ can be found in (\ref{eqn:average}).
\end{proposition}

The proof is a simple combination of the proofs of Theorem \ref{thm:nonlocal_limit} and Proposition \ref{thm:limit_myxo}.\\

\begin{figure}[h!]
\centering
\includegraphics[width=0.5\textwidth]{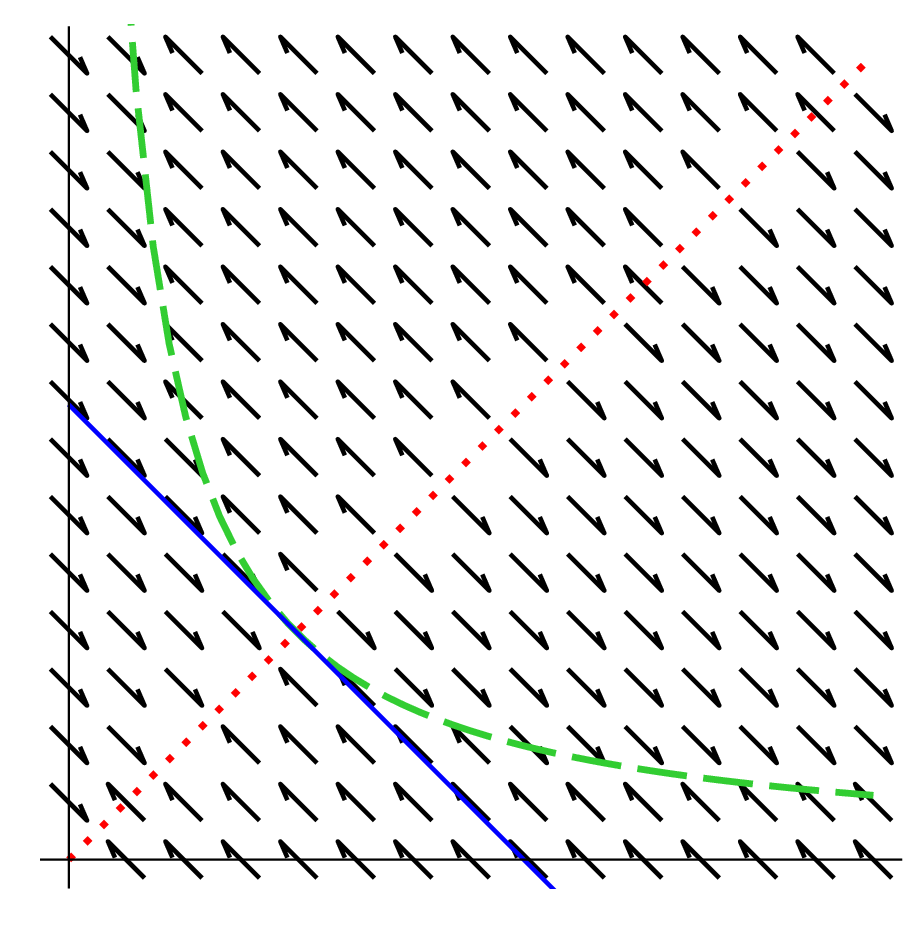}
\caption{Local dynamics for $\lambda(\rho)$ given by \eqref{eqn:lambda_ex}. The arrows mark the flow field in the $(\rho_+,\rho_-)$-plane. The red-dotted and green-dashed lines show the values for which $\lambda(\rho_+)\rho_--\lambda(\rho_-)\rho_+=0$. The blue-solid line shows the threshold values $\rho_++\rho_-=2\sqrt{\frac{\lambda_0}{\lambda_1}}$.}
\label{fig:local_dynamics}
\end{figure}

One can observe that the additional reversal term leads to a reaction term on the right-hand side of equations \eqref{eqn:macro_rho+_myxo} and \eqref{eqn:macro_rho-_myxo}, with the precise shape depending on the choice of $\lambda(\rho)$. Motivated by biological findings we expect $\lambda(0)=\lambda_0>0$ to account for spontaneous reversals in absence of other bacteria. Since head-to-head contacts increase the reversal frequency, $\lambda(\rho)$ should be an increasing function of $\rho$. To demonstrate what type of dynamics could be expected, we give an example, where we choose
\begin{align}
\label{eqn:lambda_ex}
\lambda(\rho)=\lambda_1\rho^2+\lambda_0,
\end{align}
with $\lambda_0, \lambda_1>0$. For a constant direction, $\bar \theta(x,t)\equiv \bar \theta$ and moving along the characteristics, the corresponding phase portrait is depicted in Figure \ref{fig:local_dynamics}. One can observe that wherever the total density $\rho_++\rho_->2\sqrt{\frac{\lambda_0}{\lambda_1}}$ one of the two groups dominates the other one, whilst for smaller values of $\rho_++\rho_-$ both attain the same number. This is similar to the waves observed in myxobacteria, where at the high density peaks of the waves, almost all bacteria move in the same directions, whilst in between the waves, there can be small numbers of bacteria moving in different directions.

A detailed analytical and numerical analysis for this and other choices of $\lambda(\rho)$ will be examined in a forthcoming work.

\section{Conclusion}
	
In this work we have derived macroscopic equations for self-propelled, rod shaped particles moving in two space dimensions starting at the particle level. The type of interaction considered is nematic alignment of polar particles, which is commonly found in biological and physical systems, but systematically derived hydrodynamic equations for this type of alignment are largely missing from the literature.  In order to take the limit in the hydrodynamic scaling, we used the concept of generalized collision invariants. Two types of scalings have been examined: a local scaling and a weakly non-local scaling. The resulting macroscopic equations share some similarities with their polar alignment counterparts, but differ significantly in other aspects. A key difference is the splitting of the particles into two groups moving into opposite directions. Lastly, including density dependent reversals, we have demonstrated how the equations can be used to describe the behavior of myxobacteria. The research in the area of pattern formation and non-equilibrium states of large assembles of particles is very active and more macroscopic models are needed to better understand the phase transitions and underlying causes. Several extensions of this work are possible: to describe a wider range of phenomena it could be generalized to three and more space dimensions. Further at this point the limits examined are on a formal level and lack a rigorous treatment. Extensive numerical experiments need to be conducted to compare the particle model with the macroscopic models and to characterize the different phases and states. To better understand the rippling behavior of myxobacteria different choices of reversal rates can be examined, both numerically and analytically.

\section{Appendix: proofs}
\label{sec:appendix}

\subsection{Proof of Theorem \ref{thm:limit}}
\label{subsec:app:deriv}

\subsubsection{Formal limit $\varepsilon \to 0$}
\label{subsubsec:formallimit}

All limits taken are formal. We suppose that $f^\e$ converges to some $f^0$ when $\e\rightarrow 0$. From \eqref{eqn:scaling_hydro} we see that the limiting distribution has to fulfill $Q_\text{al}f^0=0$ and therefore, use of Lemma \ref{lem:vonmises} leads us to 
\begin{align}
f^0(x,\theta,t)&=\bar f_{\rho_+(x,t), \rho_-(x,t), \bar \theta(x,t)}(\theta).
\label{eq:app:equi}
\end{align}
To derive the equations determining $\rho_+(x,t)$, $\rho_-(x,t)$, $\bar \theta(x,t)$, we multiply  \eqref{eqn:scaling_hydro} by the three GCI $\chi_{\bar \theta_{f^\varepsilon}}^\pm$ and $g_{\bar \theta_{f^\varepsilon}}$ associated to $\bar \theta_{f^\varepsilon}$ given by Prop. \ref{prop:GCI} and integrate with respect to $\theta$. 

\subsubsection{Use of GCI $\chi_{\bar \theta_{f^\varepsilon}}^\pm$}
\label{subsubsec_chipm}

We interpret the integral of $Q^0_\text{al} f^\varepsilon$ against $\chi_{\bar \theta_{f^\varepsilon}}^\pm$ as follows: 
$$\int Q^0_\text{al} f^\varepsilon (\theta) \, \chi_{\bar \theta_{f^\varepsilon}}^\pm(\theta) \dd \theta = \langle {\mathcal L}_{\bar \theta_{f^\varepsilon}} \frac{f^\varepsilon}{M_{\bar \theta_{f^\varepsilon}}} , \chi_{\bar \theta_{f^\varepsilon}}^\pm(\theta) \rangle_{{\mathcal V}_{\bar \theta_{f^\varepsilon}}',{\mathcal V}_{\theta_{f^\varepsilon}}}, $$
and, by Definition \ref{def:GCI}, this integral equals zero. Upon division by $\e$ we are left with
\begin{align*}
& \int_{ \pm \cos (\theta - \bar \theta_{f^\varepsilon}) >0 }  \big( \partial_t f^\e+v(\theta) \cdot \nabla_x f^\e \big) \dd \theta=0. 
\end{align*}
Letting $\e\rightarrow 0$, using (\ref{eq:app:equi}) and the fact that $\bar \theta_{f^\varepsilon} \to \bar \theta$, we are led to
\begin{align}
\label{eqn:left}
& \int_{ \pm \cos (\theta - \bar \theta) >0 } \big(\partial_t+v(\theta)\cdot \nabla_x \big)\, \bar f_{\rho_+(x,t), \rho_-(x,t), \bar \theta(x,t)}(\theta) \dd \theta
=0.
\end{align}
Note that the integration domain depends on $x$ and $t$.
We calculate
\begin{eqnarray*}
&&\hspace{-1cm}
\big(\partial_t + v(\theta)\cdot \nabla_x \big)\bar f_{\rho_+, \rho_-, \bar \theta}(\theta) = M_{\bar \theta} (\theta) \, \Big[ \big(\partial_t + v(\theta)\cdot \nabla_x \big) \rho_{\pm} \\
&&\hspace{0cm}
\pm \frac{\kappa \, \rho_\pm}{\cos^2(\theta - \bar \theta)} \, \sin(\theta - \bar \theta) \, \big(\partial_t + v(\theta)\cdot \nabla_x \big) \bar \theta \Big] \, \, \mbox{ for } \pm \cos(\theta - \bar \theta) >0  , 
\end{eqnarray*}
which can also be expressed as:
\begin{eqnarray}
&&\hspace{-1.5cm}
\big(\partial_t + v(\theta)\cdot \nabla_x \big)\bar f_{\rho_+, \rho_-, \bar \theta}(\theta) = M_{\bar \theta}(\theta) \Big[ \chi_{\bar \theta}^+ \big(\partial_t + v(\theta)\cdot \nabla_x \big) \rho_+  \nonumber \\
&&\hspace{-0.5cm}
+ \chi_{\bar \theta}^- \big(\partial_t + v(\theta)\cdot \nabla_x \big) \rho_-  + \kappa\frac{\sin(\theta - \bar \theta)}{\cos^2(\theta - \bar \theta)} \, \big( \chi_{\bar \theta}^+ \rho_+ - \chi_{\bar \theta}^- \rho_- \big) \, \, \big(\partial_t + v(\theta)\cdot \nabla_x \big) \bar \theta \Big] . \label{eqn:calc}
\end{eqnarray}
We note that in spite of the presence of indicator functions in (\ref{eq:def_fr+r-_2}), $\bar f_{\rho_+, \rho_-, \bar \theta}$ is continuous (and even $C^\infty$) across the points $\theta$ where $\cos(\theta - \bar \theta) =0$, so that there is no delta distribution appearing at these points when we differentiate it. In the following we perform the calculations for the integral over $\theta$ such that $\cos (\theta - \bar \theta) >0$, the case of $\cos (\theta - \bar \theta) <0$ being treated in the same manner. The left-hand side of \eqref{eqn:left}, denoted by $L$, is then equal to
\begin{eqnarray}
&&\hspace{0cm}
L= \Big( \int_{\cos (\theta - \bar \theta)>0} M_{\bar \theta} (\theta) \dd \theta \Big) \, \partial_t \rho_+ \nonumber \\
&&\hspace{0.45cm}
+ \, \Big(\int_{\cos (\theta - \bar \theta)>0} M_{\bar \theta} (\theta) \, v(\theta)  \dd \theta \Big) \cdot \nabla_x \rho_+ \nonumber \\
&&\hspace{0.45cm} + \, \kappa \rho_+  \Big( \int_{\cos (\theta - \bar \theta)>0} M_{\bar \theta} (\theta) \, \frac{\sin (\theta - \bar \theta)}{\cos^2 (\theta - \bar \theta)} \dd \theta\Big)  \,  \partial_t \bar \theta  \nonumber \\
&&\hspace{0.45cm} + \, \kappa \rho_+  \Big( \int_{\cos (\theta - \bar \theta)>0} M_{\bar \theta} (\theta) \, \frac{\sin (\theta - \bar \theta)}{\cos^2 (\theta - \bar \theta)} \, v(\theta) \dd \theta \Big) \cdot \nabla_x \bar \theta 
 \label{eqn:calc2} \\
&&\hspace{0.45cm}
= \, L_1 + \ldots + L_4 . \nonumber 
\end{eqnarray}
Thanks to the normalization condition (\ref{eq:GVM_norm}), we have $L_1 = \partial_t \rho_+$. Using that 
\begin{eqnarray}
&&
v(\theta) = \cos (\theta - \bar \theta) \, v(\bar \theta) + \sin (\theta - \bar \theta) \, v^\bot(\bar \theta), \label{eq:rotation}
\end{eqnarray}
together with cancellations due to symmetries, $L_2$ simplifies into 
\begin{align}
L_2&=\int_{\cos (\theta - \bar \theta)>0} M_{\bar \theta} (\theta) \, \big(\cos(\theta-\bar \theta) \, v(\bar \theta) \cdot \nabla_x + \sin(\theta- \bar \theta) v^\bot(\bar \theta) \cdot \nabla_x \big) \rho_+ \dd \theta\nonumber\\
&=\Big( \int_{\cos (\theta - \bar \theta)>0} M_{\bar \theta} (\theta) \, \cos(\theta - \bar \theta) \dd \theta \Big) \,  \big(v(\bar \theta) \cdot \nabla_x\big) \rho_+ \nonumber \\
&=\Big( \int_{\cos \theta >0} M (\theta) \, \cos \theta \dd \theta \Big) \,  \big(v(\bar \theta) \cdot \nabla_x\big) \rho_+ \nonumber \\
&=d_1 \,  \big(v(\bar \theta) \cdot \nabla_x\big) \rho_+,
\label{eqn:Psi01}
\end{align}
where $d_1$ is defined in (\ref{eq:def_coef}). 
We easily get $L_3 = 0$ by symmetry. 
Finally, using again (\ref{eq:rotation}), cancellations due to symmetry and integration by parts, we find that
\begin{align}
L_4&=\, \kappa \rho_+ \, \Big( \int_{\cos (\theta - \bar \theta)>0} M_{\bar \theta} (\theta) \, \frac{\sin^2 (\theta - \bar \theta)}{\cos^2 (\theta - \bar \theta)}  \dd \theta \Big) \, \big( v^\bot(\bar \theta) \cdot \nabla_x \big) \bar \theta \nonumber \\
&= - \rho_+ \, \Big( \int_{\cos (\theta - \bar \theta)>0} \frac{dM_{\bar \theta} (\theta)}{d \theta} \, \sin (\theta - \bar \theta)  \dd \theta \Big) \, \big( v^\bot(\bar \theta) \cdot \nabla_x \big) \bar \theta \nonumber \\
&= \rho_+ \, \Big( \int_{\cos (\theta - \bar \theta)>0} M_{\bar \theta} (\theta) \, \cos (\theta - \bar \theta)  \dd \theta \Big) \, \big( v^\bot(\bar \theta) \cdot \nabla_x \big) \bar \theta \nonumber \\
&= d_1 \, \rho_+ \, \big( v^\bot(\bar \theta) \cdot \nabla_x \big) \bar \theta. 
\label{eqn:Psi02}
\end{align}
Finally, collecting the contributions from above, we get (\ref{eqn:macro_rho+}) and (\ref{eqn:macro_rho-}).

\subsubsection{Use of GCI $g_{\bar \theta_{f^\varepsilon}}$}
\label{subsubsec_limpsi}

We start as in the previous section. We interpret the integral of $Q^0_\text{al} f^\varepsilon$ against $g_{\bar \theta_{f^\varepsilon}}$ as  
$$\int Q^0_\text{al} f^\varepsilon (\theta) \, g_{\bar \theta_{f^\varepsilon}} (\theta) \dd \theta =\Big \langle {\mathcal L}_{\bar \theta_{f^\varepsilon}} \frac{f^\varepsilon}{M_{\bar \theta_{f^\varepsilon}}}, g_{\bar \theta_{f^\varepsilon}}(\theta) \Big\rangle_{{\mathcal V}_{\bar \theta_{f^\varepsilon}}',{\mathcal V}_{\theta_{f^\varepsilon}}}, $$
and again by Definition \ref{def:GCI} this integral vanishes. Upon division by $\e$ we get
\begin{align*}
& \int  \big( \partial_t f^\e+v(\theta) \cdot \nabla_x f^\e \big) \, g_{\bar \theta_{f^\varepsilon}} \dd \theta=0. 
\end{align*}
Again, we assume that $f^\e \rightarrow \bar f_{\rho_+(x,t), \rho_-(x,t), \bar \theta(x,t)}$ and $\bar \theta_{f^\varepsilon} \rightarrow \bar \theta$ and find
\begin{align*}
0=&\int \Big[\big(\partial_t+v(\theta) \cdot \nabla_x\big) \bar f_{\rho_+(x,t), \rho_-(x,t), \bar \theta(x,t)}
 \Big] \, g_{\bar \theta} \dd \theta=X.
\end{align*}
We use the calculations in \eqref{eqn:calc} and split $X$ into $X=X_1 + \ldots + X_3$ with
\begin{eqnarray*}
&&\hspace{-1cm}
X_1=\int \Big[ M_{\bar \theta(x,t)} \, \big( \chi^+_{\bar \theta} \, \partial_t \rho_+ +\chi^-_{\bar \theta} \, \partial_t \rho_- \big)  \Big] \,  g_{\bar \theta}  \dd\theta, \\
&&\hspace{-1cm}
X_2=\int  M_{\bar \theta(x,t)} \, \Big[ v(\theta) \cdot \big( \chi^+_{\bar \theta} \, \nabla_x \rho_+ + \chi^-_{\bar \theta} \,\nabla_x \rho_- \big) 
+	\\
&&\hspace{6cm}
\kappa \frac{\sin (\theta - \bar \theta)}{\cos^2 (\theta - \bar \theta)} \,
\big( \chi^+_{\bar \theta} \, \rho_+ - \chi^-_{\bar \theta} \, \rho_- \big) 
\partial_t \bar \theta \Big] \, g_{\bar \theta} \dd \theta, \\
&&\hspace{-1cm}
X_3=\kappa \int M_{\bar \theta(x,t)} \, \frac{\sin (\theta - \bar \theta)}{\cos^2 (\theta - \bar \theta)} \,\big( \chi^+_{\bar \theta} \, \rho_+ - \chi^-_{\bar \theta} \, \rho_- \big) \, \big( v(\theta) \cdot \nabla_x \big) \bar \theta \, \, g_{\bar \theta} \, \dd \theta.
\end{eqnarray*}
For the first integral inside the expression of $X_1$ we have
\begin{align*}
\int M_{\bar \theta(x,t)} \, \chi^+_{\bar \theta} \,  g_{\bar \theta}  \dd\theta =0,
\end{align*}
using the symmetry relations (\ref{eq:symg}) and therefore $X_1 = 0$. 
For $X_2$, using (\ref{eq:rotation}) and cancellations due to symmetry, we get:
\begin{eqnarray*}
&&\hspace{-1cm}
X_2=\int  M_{\bar \theta(x,t)} \, \Big[ \sin(\theta - \bar \theta) \, v^\bot(\bar \theta)  \cdot \big( \chi^+_{\bar \theta} \, \nabla_x \rho_+ + \chi^-_{\bar \theta} \,\nabla_x \rho_- \big) 
+	\\
&&\hspace{6cm}
\kappa \frac{\sin (\theta - \bar \theta)}{\cos^2 (\theta - \bar \theta)} \,
\big( \chi^+_{\bar \theta} \, \rho_+ - \chi^-_{\bar \theta} \, \rho_- \big) \, 
\partial_t \bar \theta \Big] \, g_{\bar \theta} \dd \theta .
\end{eqnarray*}
Now, it turns out that the integrals with respect to $\theta$ on the two sets $\cos (\theta - \bar \theta) >0$ and $\cos (\theta - \bar \theta) <0$ are opposite by symmetry and we get
\begin{eqnarray*}
&&\hspace{-1cm}
X_2=\langle g \sin \rangle_M \, v^\bot(\bar \theta)  \cdot \nabla_x (\rho_+ - \rho_-) 
+	\kappa \, \Big\langle g \frac{\sin}{\cos^2} \Big\rangle_M \, (\rho_+ + \rho_- ) \, 
\partial_t \bar \theta .
\end{eqnarray*}
Finally for $X_3$ we again use (\ref{eq:rotation}) and cancellations due to symmetry and get:
\begin{eqnarray*}
&&\hspace{-1cm}
X_3=\kappa \int M_{\bar \theta(x,t)} \, \frac{\sin (\theta - \bar \theta)}{\cos (\theta - \bar \theta)} \,\big( \chi^+_{\bar \theta} \, \rho_+ - \chi^-_{\bar \theta} \, \rho_- \big) \, \big(  v(\bar \theta) \cdot \nabla_x \big) \bar \theta \, \, g_{\bar \theta} \, \dd \theta \\
&&\hspace{-0.35cm}
=\kappa \, \Big\langle g \frac{\sin}{\cos}\Big \rangle_M \, (\rho_+ - \rho_- ) \, 
\big(  v(\bar \theta) \cdot \nabla_x \big) \bar \theta, 
\end{eqnarray*}
since now, the integrals with respect to $\theta$ on the two sets $\cos (\theta - \bar \theta) >0$ and $\cos (\theta - \bar \theta) <0$ are equal by symmetry.
Combining these equalities, we get
\begin{align}
0&=\kappa \, \Big\langle \frac{g \sin}{\cos^2} \Big\rangle_M \, (\rho_+ + \rho_-) \, \partial_t \bar \theta + \kappa \,\Big \langle \frac{g \sin}{\cos}\Big\rangle_M \, (\rho_+ - \rho_-) \, \big( v (\bar \theta)\cdot \nabla_x \big) \bar \theta \nonumber\\
&\hspace{8cm}+ \langle g \sin \rangle_M \, v^\bot(\bar \theta)  \cdot \nabla_x (\rho_+ - \rho_-) .\label{eqn:X2X3}
\end{align}
Dividing this equation by $\kappa \,\big \langle \frac{g \sin}{\cos^2} \big\rangle_M$, we get (\ref{eqn:macro_theta}) with the coefficients defined by (\ref{eq:def_coef}).
\endproof

\subsection{Proof of Proposition \ref{thm:limit_myxo}}
\label{ssec:proof_thm:limit_myxo}

{\bf Proof.}
The proof consists of some straightforward extensions to the proof of Theorem \ref{thm:limit}. 
First of all we note that the limiting distribution $f^0$ still has to fulfill $Q^0_\text{al}f^0=0$. Therefore we get that
\begin{align*}
f^0(x,\theta, t)=\bar f_{\rho_+(x,t), \rho_-(x,t), \bar \theta(x,t)}(\theta)
\end{align*}
and we are left with determining the integral of \eqref{eqn:Qrev} against the generalized collision invariants $\chi_{\bar \theta_{f^\varepsilon}}^\pm$ and $g_{\bar \theta_{f^\varepsilon}}$. We start by noting that
\begin{eqnarray}
Q^0_\text{rev} \bar f_{\rho_+, \rho_-, \bar \theta} (\theta) =\begin{cases}
\big( \lambda(\rho_+) \, \rho_- - \lambda(\rho_-) \, \rho_+ \big) \, M_{\bar \theta}(\theta) \, \,  \, \, \text{for} \,  \cos(\theta-\bar \theta) > 0, 
\\
\big( \lambda(\rho_-) \, \rho_+ -\lambda(\rho_+) \, \rho_- \big) \, M_{\bar \theta}(\theta) \, \,  \, \, \text{for} \,  \cos(\theta-\bar \theta) < 0. 
\end{cases} 
\label{eqn:Qrevfbar}
\end{eqnarray}
With this it follows easily that 
\begin{align}
\label{eqn:left_mxyo}
& \int Q_\text{rev}^0 \bar f_{\rho_+(x,t), \rho_-(x,t), \bar \theta(x,t)}\chi_{\bar \theta_{f^\varepsilon}}^\pm \dd \theta
= \pm \big(\lambda(\rho_+)\rho_--\lambda(\rho_-)\rho_+\big).
\end{align}
Finally following the same arguments as in Section \ref{subsubsec_limpsi} about $X_1$, we find that
\begin{align*}
\int Q^0_\text{rev} \bar f_{\rho_+(x,t), \rho_-(x,t), \bar \theta(x,t)} \,  \, g_{\bar \theta}  \, \dd\theta  =0,
\end{align*}
which finishes the proof.
\endproof

\subsection{Weakly non-local scaling: proof of Theorem \ref{thm:nonlocal_limit}}
\label{app:nonlocal}

Compared to the proof of Theorem \ref{thm:limit}, we need to compute the contributions of 
\begin{eqnarray}
&&\hspace{-1cm}
T_1^\pm = \int Q^1_\text{al} \bar f_{\rho_+, \rho_-, \bar \theta}(\theta) \, \chi_{\bar \theta}^{\pm} (\theta)  \dd \theta, \quad 
T_2 = \int Q^1_\text{al} \bar f_{\rho_+, \rho_-, \bar \theta}(\theta) \, g_{\bar \theta} (\theta)  \dd \theta. \label{def_T1T2}
\end{eqnarray}
We first compute:
\begin{eqnarray*}
j_{\bar f_{\rho_+, \rho_-, \bar \theta}} &=& \rho_+ \int_{\cos(\theta - \bar \theta) >0} v(2 \theta) M_{\bar \theta}(\theta)  \dd \theta + 
\rho_- \int_{\cos(\theta - \bar \theta) <0} v(2 \theta) M_{\bar \theta}(\theta)  \dd \theta. 
\end{eqnarray*}
Thanks to 
$$ v(2\theta) = \cos \big(2 (\theta - \bar \theta)\big) \, v(2 \bar \theta) + \sin \big(2 (\theta - \bar \theta)\big) \, v^\bot(2 \bar \theta), $$
and owing to cancellations by symmetry, we get
\begin{eqnarray*}
j_{\bar f_{\rho_+, \rho_-, \bar \theta}} &=& (\rho_+ + \rho_-) \, \left( \int_{\cos \theta >0}  \cos (2 \theta) \, M(\theta) \, \dd \theta \right) \, v(2 \bar \theta)\\
&=& (\rho_+ + \rho_-) \, \langle \cos (2 \cdot) \rangle_M \,  v(2 \bar \theta), 
\end{eqnarray*}
and consequently
\begin{eqnarray*}
\big| j_{\bar f_{\rho_+, \rho_-, \bar \theta}} \big|&=& (\rho_+ + \rho_-) \, \langle \cos (2 \cdot) \rangle_M ,
\end{eqnarray*}
since $\langle \cos (2 \cdot) \rangle_M >0$. Now we get:
\begin{eqnarray*}
j^1_{\bar f_{\rho_+, \rho_-, \bar \theta}} &=& k \, \Delta_x j_{\bar f_{\rho_+, \rho_-, \bar \theta}} \\
&=& k \, \langle \cos (2 \cdot) \rangle_M \, \Delta_x \big( (\rho_+ + \rho_-) \,  v(2 \bar \theta) \big) \\
&=&  k \, \langle \cos (2 \cdot) \rangle_M \, \Big( \Delta_x  (\rho_+ + \rho_-) \,  v(2 \bar \theta) + 4 \nabla_x  (\rho_+ + \rho_-) \cdot \nabla_x \bar \theta \, \, v^\bot (2 \bar \theta)   \\
&&\hspace{4cm}
+4 (\rho_+ + \rho_-) \,  \big( \Delta_x \bar \theta \, v^\bot (2 \bar \theta) - |\nabla_x \bar \theta |^2 \,  v(2 \bar \theta) \big) \Big), 
\end{eqnarray*}
and finally
\begin{eqnarray}
\bar \theta_{\bar f_{\rho_+, \rho_-, \bar \theta}}^1&=& \frac{1}{2|j_{\bar f_{\rho_+, \rho_-, \bar \theta}}|} \big(j_{\bar f_{\rho_+, \rho_-, \bar \theta}}^1 \cdot v^\bot(2 \bar \theta_{\bar f_{\rho_+, \rho_-, \bar \theta}}) \big) \nonumber \\
&=& \frac{2k}{\rho_+ + \rho_-} \, \big( \nabla_x (\rho_+ + \rho_-) \cdot \nabla_x \bar \theta +  (\rho_+ + \rho_-) \, \Delta_x \bar \theta \big)
\nonumber \\
&=& \frac{2k}{\rho_+ + \rho_-} \,  \nabla_x \cdot \big( (\rho_+ + \rho_-)\nabla_x \bar \theta \big) . \label{eq:bar_theta_1}
\end{eqnarray}
Now we have 
\begin{eqnarray*}
&&\hspace{-1cm}
Q^1_\text{al} \bar f_{\rho_+, \rho_-, \bar \theta}(\theta)  \nonumber \\
&&\hspace{-0.5cm}
= \bar \theta_{\bar f_{\rho_+, \rho_-, \bar \theta}}^1 \, \partial_\theta  \bigg[ \Big( 2 \, \delta \big(\cos(\theta- \bar \theta) \big) \, \sin^2 (\theta- \bar \theta) 
- \big|\cos(\theta- \bar \theta) \big| \, \Big)  \, M_{\bar \theta} \big( \rho_+ \chi_{\bar \theta}^+ + \rho_- \chi_{\bar \theta}^- \big)\nonumber \\
&&\hspace{3cm}
+ 2D \,  \cos(\theta- \bar \theta) \, \sin(\theta- \bar \theta) \, \partial_\theta \Big( M_{\bar \theta} \big( \rho_+ \chi_{\bar \theta}^+ + \rho_- \chi_{\bar \theta}^- \big) \Big) \bigg],
\end{eqnarray*}
however, 
\begin{equation}
\delta \big(\cos(\theta- \bar \theta) \big) \, M_{\bar \theta} = 0, 
\label{eq:deltaM}
\end{equation} 
and similarly 
$$ \partial_\theta \Big( M_{\bar \theta} \big( \rho_+ \chi_{\bar \theta}^+ + \rho_- \chi_{\bar \theta}^- \big) \Big)
= \big( \rho_+ \chi_{\bar \theta}^+ + \rho_- \chi_{\bar \theta}^- \big) \, \partial_\theta   M_{\bar \theta}. \, 
$$
Together with 
$$ \partial_\theta M_{\bar \theta} (\theta) = - \kappa \frac{\sin(\theta - \bar\theta)}{\cos^2(\theta - \bar\theta)} \, \mbox{Sign} \big( \cos(\theta - \bar\theta) \big) M_{\bar \theta} (\theta) , $$
we find (noting that $\kappa D = 1$):
\begin{eqnarray*}
&&\hspace{-1cm}
Q^1_\text{al} \bar f_{\rho_+, \rho_-, \bar \theta}(\theta)  
= - \bar \theta_{\bar f_{\rho_+, \rho_-, \bar \theta}}^1 \, \partial_\theta  \left[ 
\bigg( \big|\cos(\theta- \bar \theta) \big|  + 2 \frac{\sin^2(\theta- \bar \theta)}{\big|\cos(\theta- \bar \theta) \big|} \bigg)  \, \big( \rho_+ \chi_{\bar \theta}^+ + \rho_- \chi_{\bar \theta}^- \big) \, M_{\bar \theta} \right].
\end{eqnarray*}

We first compute $T_1^+$ (the computation is obviously similar for $T_1^-$). We have 
\begin{eqnarray*}
&&\hspace{-1cm}
T_1^+ = \int_{\cos(\theta-\bar \theta) >0} Q^1_\text{al} \bar f_{\rho_+, \rho_-, \bar \theta}(\theta) \dd \theta  \\
&&\hspace{-1cm}
= - \bar \theta_{\bar f_{\rho_+, \rho_-, \bar \theta}}^1 \, \int_{\cos(\theta-\bar \theta) >0} \partial_\theta  \left[ 
\bigg( \big|\cos(\theta- \bar \theta) \big|  + 2 \frac{\sin^2(\theta- \bar \theta)}{\big|\cos(\theta- \bar \theta) \big|} \bigg)  \, \big( \rho_+ \chi_{\bar \theta}^+ + \rho_- \chi_{\bar \theta}^- \big) \, M_{\bar \theta} \right]  \dd\theta \\
&&\hspace{-1cm}
= - \bar \theta_{\bar f_{\rho_+, \rho_-, \bar \theta}}^1 \,   \bigg[ 
\bigg( \big|\cos(\theta- \bar \theta) \big|  + 2 \frac{\sin^2(\theta- \bar \theta)}{\big|\cos(\theta- \bar \theta) \big|} \bigg)  \, \big( \rho_+ \chi_{\bar \theta}^+ + \rho_- \chi_{\bar \theta}^- \big) \, M_{\bar \theta} \bigg]_{\theta - \bar \theta = - \pi/2}^{\theta - \bar \theta = + \pi/2} = 0, 
\end{eqnarray*}
because the contact of $M_{\bar \theta}$ to $0$ at $\theta - \bar \theta = \pm \pi/2$ is infinite order and compensates for the blow up $|\cos(\theta -  \bar \theta)|$ at the denominator. 

We are left with the computation of $T_2$. Using (\ref{eq:transinvar}) and (\ref{eq:deltaM}), we can write:
\begin{align*}
T_2& = \int_{-\pi}^\pi Q^1_\text{al} \bar f_{\rho_+, \rho_-, \bar \theta}(\theta)  \, g_{\bar \theta}(\theta) \dd  \theta  \\
&= - \bar \theta_{\bar f_{\rho_+, \rho_-, \bar \theta}}^1 \, \int_{-\pi}^\pi \partial_\theta  \left[ 
\left( \big|\cos(\theta- \bar \theta) \big|  + 2 \frac{\sin^2(\theta- \bar \theta)}{\big|\cos(\theta- \bar \theta) \big|} \right)  \, \big( \rho_+ \chi_{\bar \theta}^+ + \rho_- \chi_{\bar \theta}^- \big) \, M_{\bar \theta} \right]  \, g_{\bar \theta}(\theta) \dd\theta \\
&= - \bar \theta_{\bar f_{\rho_+, \rho_-, \bar \theta}}^1 \, \int_{-\pi}^\pi \partial_\theta  \left[ 
\left( |\cos \theta |  + 2 \frac{\sin^2 \theta}{|\cos \theta |} \right)  \, \big( \rho_+ \chi^+ + \rho_- \chi^- \big) \, M \right]  \, g(\theta) \dd\theta \\
&= - \bar \theta_{\bar f_{\rho_+, \rho_-, \bar \theta}}^1 \, \int_{-\pi}^\pi \partial_\theta  \left[ 
\left( |\cos \theta |  + 2 \frac{\sin^2 \theta}{|\cos \theta |} \right)  \, M \right]  \, g(\theta) \, \big( \rho_+ \chi^+ + \rho_- \chi^- \big) \dd\theta. 
\end{align*}
We have 
\begin{eqnarray*}
&&\hspace{-1cm}
\partial_\theta  \left[ \left( |\cos \theta |  + 2 \frac{\sin^2 \theta}{|\cos \theta |} \right)  \, M  \, \right] = \bigg\{ \mbox{Sign} \big( \cos \theta \big)  \Big( \sin \theta   + 2 \frac{\sin \theta}{\cos^2 \theta } \Big) - \kappa \left( \frac{\sin \theta}{\cos \theta } +2 \frac{\sin^3 \theta}{\cos^3 \theta } \right) \bigg\} \, M . 
\end{eqnarray*}
We note that this expression is changed into its opposite under the change $\theta \to \pi - \theta$. Since $g \, M$ is also changed into its opposite under this transformation, it follows that the factors of $\rho_+$ and $\rho_-$ are equal, which leads to:
\begin{eqnarray*}
&&\hspace{-1cm}
T_2 = - \bar \theta_{\bar f_{\rho_+, \rho_-, \bar \theta}}^1 \,  
\bigg[ \langle \sin \theta \, g \rangle_M + 2 \Big\langle \frac{\sin \theta}{\cos^2 \theta} \, g \Big\rangle_M  \\
&&\hspace{5cm}
- \kappa \,  \left( \Big\langle \frac{\sin \theta}{\cos \theta} \, g \Big\rangle_M + 2 \Big\langle \frac{\sin^3 \theta}{\cos^3 \theta} \, g \Big\rangle_M \right) \bigg]  \, ( \rho_+  + \rho_-) . 
\end{eqnarray*}
Finally, inserting (\ref{eq:bar_theta_1}) into this expression leads to (\ref{eqn:macro_nonlocal}), which ends the proof. \endproof


\begin{thebibliography}{10}

\bibitem{Aoki1982}
I.~Aoki.
\newblock A simulation study on the schooling mechanism in fish.
\newblock {\em Bull. Jpn. Soc. Sci. Fish.}, 48(8):1081--1088, 1982.

\bibitem{Barbaro_Degond_2014}
A.~Barbaro and P.~Degond.
\newblock Phase transition and diffusion among socially interacting
  self-propelled agents.
\newblock {\em Discrete Contin. Dyn. Syst. Ser. B}, 19:1249--1278, 2014.

\bibitem{Baskaran2008}
A.~Baskaran and M.~Marchetti.
\newblock Hydrodynamics of self-propelled hard rods.
\newblock {\em Phys. Rev. E}, 77(1):011920, 2008.

\bibitem{Baskaran2010}
A.~Baskaran and M.~Marchetti.
\newblock Nonequilibrium statistical mechanics of self-propelled hard rods.
\newblock {\em J. Stat. Mech.: Theory Exp.}, (4):1--14, 2010.

\bibitem{Ben-Jacob2000}
E.~Ben-Jacob, I.~Cohen, and H.~Levine.
\newblock Cooperative self-organization of microorganisms.
\newblock {\em Adv. in Phys.}, 49(4):395--554, 2000.

\bibitem{Bertin2009}
E.~Bertin, M.~Droz, and G.~Gr\'{e}goire.
\newblock Hydrodynamic equations for self-propelled particles: microscopic
  derivation and stability analysis.
\newblock {\em J. Phys. A: Math. Theor.}, 42(44):445001, 2009.

\bibitem{Borner2002}
U.~B\"{o}rner, A.~Deutsch, H.~Reichenbach, and M.~B\"{a}r.
\newblock Rippling patterns in aggregates of myxobacteria arise from cell-cell
  collisions.
\newblock {\em Phys. Rev. Lett.}, 89:078101, 2002.

\bibitem{Buhl2006}
J.~Buhl, D.~Sumpter, I.~Couzin, J.~Hale, E.~Despland, E.~Miller, and
  S.~Simpson.
\newblock From disorder to order in marching locusts.
\newblock {\em Science}, 312(5778):1402--1406, 2006.

\bibitem{Carrillo2009}
J.~Carrillo, M.~D\textsc{\char13}Orsogna, and V.~Panferov.
\newblock Double milling in self-propelled swarms from kinetic theory.
\newblock {\em Kinet. Relat. Models}, 2(2):363--378, 2009.

\bibitem{Cavagna2010}
A.~Cavagna, A.~Cimarelli, I.~Giardina, G.~Parisi, R.~Santagati, F.~Stefanini,
  and M.~Viale.
\newblock Scale-free correlations in starling flocks.
\newblock {\em Proc. Natl. Acad. Sci. USA}, 107(26):11865--11870, 2010.

\bibitem{Chate2008}
H.~Chat\'{e}, F.~Ginelli, G.~Gr\'{e}goire, and F.~Raynaud.
\newblock Collective motion of self-propelled particles interacting without
  cohesion.
\newblock {\em Phys. Rev. E}, 77(4):046113, 2008.

\bibitem{Chuang2007}
Y.~Chuang, M.~D\textsc{\char13}Orsogna, D.~Marthaler, A.~Bertozzi, and
  L.~Chayes.
\newblock State transitions and the continuum limit for a 2{D} interacting,
  self-propelled particle system.
\newblock {\em Phys. D}, 232(1):33--47, 2007.

\bibitem{Couzin2002}
I.~Couzin, J.~Krause, R.~James, G.~Ruxton, and N.~Franks.
\newblock Collective memory and spatial sorting in animal groups.
\newblock {\em J. Theoret. Biol.}, 218(1):1--11, 2002.

\bibitem{DeGennes_Prost_1993}
P.~G. de~Gennes and J.~Prost.
\newblock {\em The physics of liquid crystals}.
\newblock Oxford University Press, 1993.

\bibitem{Degond_Delebecque_Peurichard_2015}
P.~Degond, F.~Delebecque, and D.~Peurichard.
\newblock Continuum model for linked fibers with alignment interactions.
\newblock {\em Math. Models Methods Appl. Sci.}, to appear, 2015.

\bibitem{Degond_Dimarco_Mac_2014}
P.~Degond, G.~Dimarco, and T.~B.~N. Mac.
\newblock {Hydrodynamics of the Kuramoto-Vicsek model of rotating
  self-propelled particles}.
\newblock {\em Math. Models Methods Appl. Sci.}, 24:277--325, 2014.

\bibitem{Degond_etal_2015}
P.~Degond, G.~Dimarco, T.~B.~N. Mac, and N.~Wang.
\newblock Macroscopic models of collective motion with repulsion.
\newblock {\em Commun. Math. Sci}, 13:1615--1638, 2015.

\bibitem{Degond_Frouvelle_Liu_2013}
P.~Degond, A.~Frouvelle, and J.-G. Liu.
\newblock Macroscopic limits and phase transition in a system of self-propelled
  particles.
\newblock {\em J. Nonlinear Sci.}, 23:427--456, 2013.

\bibitem{Degond_Frouvelle_Liu_2015}
P.~Degond, A.~Frouvelle, and J.-G. Liu.
\newblock Phase transitions, hysteresis, and hyperbolicity for self-organized
  alignment dynamics.
\newblock {\em Arch. Ration. Mech. Anal}, 216:63--115, 2015.

\bibitem{Degond_Liu_2012}
P.~Degond and J.-G. Liu.
\newblock {Hydrodynamics of self-alignment interactions with precession and
  derivation of the Landau-Lifschitz-Gilbert equation}.
\newblock {\em Math. Models Methods Appl. Sci.}, 22 Suppl. 1:114001, 2012.

\bibitem{Degond2013}
P.~Degond, J.-G. Liu, S.~Motsch, and V.~Panferov.
\newblock Hydrodynamic models of self-organized dynamics: derivation and
  existence theory.
\newblock {\em Methods Appl. Anal.}, 20(2):89--114, 2013.

\bibitem{Degond_Liu_Ringhofer_2014}
P.~Degond, J.-G. Liu, and C.~Ringhofer.
\newblock Evolution of wealth in a nonconservative economy driven by local nash
  equilibria.
\newblock {\em Philos. Trans. A}, 372:20130394, 2015.

\bibitem{Degond2008}
P.~Degond and S.~Motsch.
\newblock Continuum limit of self-driven particles with orientation
  interaction.
\newblock {\em Math. Models Methods Appl. Sci.}, 18(Suppl.):1193--1215, 2008.

\bibitem{Degond_Navoret_2015}
P.~Degond and L.~Navoret.
\newblock A multi-layer model for self-propelled disks interacting through
  alignment and volume exclusion.
\newblock {\em Math. Models Methods Appl. Sci.}, to appear, 2015.

\bibitem{Degond_Hui_2015}
P.~Degond and H.~Yu.
\newblock Self-organized hydrodynamics in an annular domain: modal analysis and
  nonlinear effects.
\newblock {\em Math. Models Methods Appl. Sci.}, 25:495--519, 2015.

\bibitem{Dhar2006}
P.~Dhar, T.~Fischer, Y.~Wang, T.~Mallouk, W.~Paxton, and A.~Sen.
\newblock Autonomously moving nanorods at a viscous interface.
\newblock {\em Nano Lett.}, 6(1):66--72, 2006.

\bibitem{Dworkin1996}
M.~Dworkin.
\newblock Recent advances in the social and developmental biology of the
  myxobacteria.
\newblock {\em Microbiol. Rev.}, 60(1):70--102, 1996.

\bibitem{Frouvelle2012}
A.~Frouvelle.
\newblock A continuum model for alignment of self-propelled particles with
  anisotropy and density-dependent parameters.
\newblock {\em Math. Models Methods Appl. Sci.}, 22(7):1250011, 2012.

\bibitem{Ginelli2010}
F.~Ginelli, F.~Peruani, M.~B{\"a}r, and H.~Chat{\'e}.
\newblock Large-scale collective properties of self-propelled rods.
\newblock {\em Phys. Rev. Lett.}, 104(18):184502.

\bibitem{Helbing2007}
D.~Helbing, A.~Johansson, and H.~Al-Abideen.
\newblock Dynamics of crowd disasters: an empirical study.
\newblock {\em Phys. Rev. E}, 75(4):046109, 2007.

\bibitem{Igoshin2001}
O.~Igoshin, A.~Mogilner, R.~Welch, D.~Kaiser, and G.~Oster.
\newblock Pattern formation and traveling waves in myxobacteria: theory and
  modeling.
\newblock {\em Proc. Natl. Acad. Sci. USA}, 98(26):14913--14918, 2001.

\bibitem{Igoshin2004}
O.~Igoshin, R.~Welch, D.~Kaiser, and G.~Oster.
\newblock Waves and aggregation patterns in myxobacteria.
\newblock {\em Proc. Natl. Acad. Sci. USA}, 101(12):4256--4261, 2004.

\bibitem{Kaiser2003}
D.~Kaiser.
\newblock Coupling cell movement to multicellular development in myxobacteria.
\newblock {\em Nat. Rev. Microbiol.}, 1(1):45--54, 2003.

\bibitem{Mogilner1999}
A.~Mogilner and L.~Edelstein-Keshet.
\newblock A non-local model for a swarm.
\newblock {\em J. Math. Biol.}, 38(6):534--570, 1999.

\bibitem{Moussaid2010}
M.~Moussa\"{\i}d, N.~Perozo, S.~Garnier, D.~Helbing, and G.~Theraulaz.
\newblock The walking behaviour of pedestrian social groups and its impact on
  crowd dynamics.
\newblock {\em PLoS ONE}, 5(4):1--7, 2010.

\bibitem{Jiang_Xiong_Zhang_2015}
J.~N., L.~Xiong, and T.-F. Zhang.
\newblock Hydrodynamic limits of the self-organized hydrodynamic models.
\newblock {\em arXiv}, 1508.04640, 2015.

\bibitem{Narayan2007}
V.~Narayan, S.~Ramaswamy, and N.~Menon.
\newblock Long-lived giant number fluctuations in a swarming granular nematic.
\newblock {\em Science}, 317(5834):105--108, 2007.

\bibitem{Ngo2012}
S.~Ngo, F.~Ginelli, and H.~Chat{\'e}.
\newblock Competing ferromagnetic and nematic alignment in self-propelled polar
  particles.
\newblock {\em Phys. Rev. E}, 86(5):050101, 2012.

\bibitem{Parrish1997}
J.~Parrish and W.~Hamner.
\newblock {\em Animal groups in three dimensions: how species aggregate}.
\newblock Cambridge University Press, 1997.

\bibitem{Peruani2011}
F.~Peruani, F.~Ginelli, M.~B{\"a}r, and H.~Chat{\'e}.
\newblock Polar vs. apolar alignment in systems of polar self-propelled
  particles.
\newblock {\em J. Phys. Conf. Ser.}, 297:012014, 2011.

\bibitem{Reynold1987}
C.~Reynold.
\newblock Flocks, herds, and schools: A distributed behavioral model.
\newblock {\em SIGGRAPH Comput. Graph.}, 21(4):25--34, 1987.

\bibitem{Sokolov2007}
A.~Sokolov, I.~Aranson, J.~Kessler, and R.~Goldstein.
\newblock Concentration dependence of the collective dynamics of swimming
  bacteria.
\newblock {\em Phys. Rev. Lett.}, 98(15):158102, 2007.

\bibitem{Toner_Tu_1995}
J.~Toner and Y.~Tu.
\newblock Long-range order in a two-dimensional dynamical xy model: how birds
  fly together.
\newblock {\em Phys. Rev. Lett.}, 75:4326--4329, 1995.

\bibitem{Vicsek1995}
T.~Vicsek, A.~Czir{\'o}k, E.~Ben-Jacob, I.~Cohen, and O.~Shochet.
\newblock Novel type of phase transition in a system of self-driven particles.
\newblock {\em Phys. Rev. Lett.}, 75(6):1226--1229, 1995.

\bibitem{Vicsek_Zafeiris_PhysRep12}
T.~Vicsek and A.~Zafeiris.
\newblock Collective motion.
\newblock {\em Phys. Rep.}, (517):71--140, 2012.

\bibitem{Welch2001}
R.~Welch and D.~Kaiser.
\newblock Cell behavior in traveling wave patterns of myxobacteria.
\newblock {\em Proc. Natl. Acad. Sci. USA}, 98(26):14907--14912, 2001.

\bibitem{Wu2009}
Y.~Wu, A.~Kaiser, Y.~Jiang, and M.~Alber.
\newblock Periodic reversal of direction allows myxobacteria to swarm.
\newblock {\em Proc. Natl. Acad. Sci. USA}, 106:1222--1227, 2009.

\bibitem{Zhang2010}
H.-P. Zhang, A.~Be'er, E.-L. Florin, and H.~Swinney.
\newblock Collective motion and density fluctuations in bacterial colonies.
\newblock {\em Proc. Natl. Acad. Sci. USA}, 107(31):13626--13630, 2010.

\end{thebibliography}

\end{document}